\documentclass[10pt]{article}

\usepackage{latexsym}
\usepackage{amsfonts}
\usepackage{amsmath}
\usepackage{mathrsfs}  
\usepackage{dsfont}    
\usepackage{bbold}     
\usepackage[english]{babel}
\usepackage{caption}
\usepackage{epsfig}
\usepackage{float}
\usepackage{subfigure}

\newtheorem{theorem}{Theorem}

\newtheorem{lemma}{Lemma}

\newcommand{\zerarcounters}{\setcounter{equation}{0}}

\newcommand{\ZZZ}{\mathds{Z}}

\newcommand{\NNN}{\mathds{N}}

\newcommand{\RRR}{\mathds{R}}

\newcommand{\AAA}{{\mathcal A}}
\newcommand{\BB}{{\mathcal B}}
\newcommand{\CC}{{\mathcal C}}
\newcommand{\DD}{{\mathcal D}}

\newcommand{\HHH}{{\mathcal H}}

\newcommand{\PPP}{{\mathcal P}}

\newcommand{\SSS}{{\mathcal S}}
\newcommand{\TT}{{\mathcal T}}

\newcommand{\Fullbox}{{\rule{2.0mm}{2.0mm}}}

\newcommand{\EP}{\hfill\Fullbox\vspace{0.2cm}}
\newcommand{\prova}{\noindent{\it Proof. }}
\newcommand{\io}{\infty}
\newcommand{\eps}{\varepsilon}
\newcommand{\al}{\alpha}
\newcommand{\be}{\beta}

\newcommand{\n}{\nu}

\newcommand{\x}{\xi}

\newcommand{\g}{\gamma}

\newcommand{\f}{\varphi}

\newcommand{\oo}{\boldsymbol{\omega}}
\newcommand{\nn}{\boldsymbol{\nu}}
\newcommand{\pps}{\boldsymbol{\psi}}
\newcommand{\vzero}{\boldsymbol{0}}

\def\vec#1{\mbox{\boldmath $#1$}}
\def\ins#1#2#3{\vbox to0pt{\kern-#2 \hbox{\kern#1 #3}\vss}\nointerlineskip}

\begin{document}

\title{\bf Globally and locally attractive solutions\\
for quasi-periodically forced systems}
\author
{\bf Michele V. Bartuccelli$^\ast$, Jonathan H.B. Deane$^\ast$
and Guido Gentile$^\dagger$
\vspace{2mm} \\ \small 
$^\ast$Department of Mathematics and Statistics,
University of Surrey, Guildford, GU2 7XH, UK.
\\ \small 
E-mails:
m.bartuccelli@surrey.ac.uk, j.deane@surrey.ac.uk
\\ \small 
$^\dagger$Dipartimento di Matematica, Università di Roma Tre, Roma,
I-00146, Italy.
\\ \small 
E-mail: gentile@mat.unirom3.it
}

\date{}

\maketitle

\begin{abstract}
We consider a class of differential equations, $\ddot x +
\gamma \dot x + g(x) = f(\oo t)$, with $\oo \in \RRR^{d}$,
describing one-dimensional dissipative systems subject to a periodic
or quasi-periodic (Diophantine) forcing. We study existence and
properties of the limit cycle described by the trajectory with the same
quasi-periodicity as the forcing. For $g(x)=x^{2p+1}$, $p\in\NNN$,
we show that, when the dissipation coefficient is large enough,
there is only one limit cycle and that it is a global attractor.
In the case of other forces, including $g(x)=x^{2p}$
(with $p=1$ describing the varactor equation), we find estimates
for the basin of attraction of the limit cycle.
\end{abstract}




\zerarcounters
\section{Introduction}\label{sec:1}

Consider the ordinary differential equation
\begin{equation} \ddot x + \g \dot x + x^{2p+1} = f(\oo t) ,
\label{eq:1.1} \end{equation}
where $p \in \NNN$, $\oo\in\RRR^{d}$ is the frequency vector,
$f(\pps)$ is an analytic quasi-periodic function,
\begin{equation} f(\pps) = \sum_{\nn\in \ZZZ^{d}}
{\rm e}^{i\nn \cdot \pps} f_{\nn} ,
\label{eq:1.2} \end{equation}
with average $\langle f \rangle \equiv f_{\vzero} \neq 0$,
and $\g>0$ is a real parameter (dissipation coefficient).
Here and henceforth we denote with $\cdot$ the scalar product
in $\RRR^{d}$. By the analyticity assumption on $f$ there are
two strictly positive constants $F$ and $\x$ such that
one has $|f_{\nn}|\le F \, {\rm e}^{-\x|\nn|}$ for all $\nn\in\ZZZ^{d}$.

If $d>1$ we shall assume a Diophantine condition on the
frequency vector $\oo$, that is
\begin{equation} \left| \oo\cdot\nn \right| \ge C_{0} |\nn|^{-\tau} \qquad
\forall \nn\in\ZZZ^{d} \setminus\{\vzero\} ,
\label{eq:1.3} \end{equation}
where $|\nn|=|\nn|_{1}\equiv|\n_{1}|+\ldots+|\n_{d}|$, and $C_{0}$
and $\tau$ are positive constants, with $\tau>d-1$ and $C_{0}$
small enough. Note that for $d=1$ the condition~(\ref{eq:1.3}) is
automatically satisfied for all $\oo$.

In this paper we want to show that for $\g$ large enough
the system~(\ref{eq:1.1}) admits a global attractor
which is a quasi-periodic solution with the
same frequency vector $\oo$ as the forcing $f$.
This will be done in two steps: first we prove that for $\g$ large enough
there is a quasi-periodic solution $x_{0}(t)$ with frequency vector $\oo$
(cf. Theorem~\ref{thm:1} in Section~\ref{sec:2}), second we prove that,
again for $\g$ large enough, any trajectory is attracted by $x_{0}(t)$
(cf. Theorem~\ref{thm:2} in Section~\ref{sec:3}).

In particular, this solves for the system~(\ref{eq:1.1}) a problem
left as open in~\cite{GBD}. Indeed in~\cite{GBD} we considered
a class of ordinary differential equations, including~(\ref{eq:1.1}),
and proved existence of a quasi-periodic solution with the same
quasi-periodicity as the forcing, but we couldn't conclude that
this was the only solution with such a property. The result stated
above gives an affirmative answer to this problem for
the system~(\ref{eq:1.1}), by showing that
the quasi-periodic solution $x_{0}(t)$ is unique;
cf. Theorem~\ref{thm:3} in Section~\ref{sec:4}.

This uniqueness result holds for the more general systems
studied in~\cite{GBD}, including the resistor-inductor-varactor circuit,
or simply varactor equation, described in~\cite{GBD,BDGM}.
In that case the solution $x_{0}(t)$ is not a global attractor
any more, but it turns out to be the only attractor in a
neighbourhood of the solution itself.

More precisely the situation is as follows.
We can consider systems described by
\begin{equation}
\ddot x + \g \dot x + g(x) =  f(\oo t) ,
\label{eq:1.4} \end{equation}
where $f$ is given by~(\ref{eq:1.2}) and $g$ is an analytic function.
Studying the behaviour of the system~(\ref{eq:1.4}) for $\g$ large
enough suggests to introduce a new parameter $\eps=1/\g$,
in terms of which the differential equation~(\ref{eq:1.4}) becomes
\begin{equation}
\eps \ddot x + \dot x + \eps g(x) =  \eps f(\oo t) ,
\label{eq:1.5} \end{equation}
and study what happens for $\eps$ small enough.

If we assume that there exists $c_{0}\in\RRR$ such that
$g(c_{0})=f_{\vzero}$ and $g'(c_{0}):=\partial_{x}g(c_{0})\neq 0$,
then the system~(\ref{eq:1.5}) admits a quasi-periodic solution $x_{0}(t)$,
analytic in $t$, with the same frequency vector $\oo$ as the forcing $f$,
and furthermore $x_{0}(t)=c_{0}+O(\eps)$. This was proved in~\cite{GBD},
where the solution $x_{0}(t)$ was explicitly constructed through a
suitable summation of the perturbation series
\begin{equation}
x_{P} (t) = \sum_{k=0}^{\io} \eps^{k} x^{(k)}(\oo t) ,
\qquad x^{(k)}(\pps) =
\sum_{\nn\in\ZZZ^{d}} {\rm e}^{i\nn\cdot\pps} x^{(k)}_{\nn} ,
\label{eq:1.6} \end{equation}
for a function formally solving the equations of motion.

As a drawback of the construction we were not able to prove any
uniqueness result about $x_{0}(t)$. In fact, in principle, there could
be other quasi-periodic solutions near $x_{0}(t)$, possibly with the
same frequency vector $\oo$. Neither we could exclude existence of
other solutions reducing to $c_{0}$ as $\eps\to0$ or even admitting
the same formal expansion~(\ref{eq:1.6}) in powers of $\eps$.
In this papers we get rid of these possibilities, and we prove that
there exists, in the plane $(x,\dot x)$, a neighborhood $\BB$ of the
point $(c_{0},0)$ where $(x_{0}(t),\dot x_{0}(t))$ is the only stable
solution of~(\ref{eq:1.5}). Moreover it turns out to be asymptotically
stable, that is it attracts any trajectory starting in $\BB$.
Therefore, this allows us to formulate a strengthened
version of the theorem proved in~\cite{GBD};
cf. Theorem~\ref{thm:4} in Section~\ref{sec:4}.

In general the neighbourhood $\BB$ can be very small.
In specific cases one can look for improved estimates of $\BB$.
In particular we can consider the case $g(x)=x^{2}$, corresponding to
the varactor equation studied in~\cite{BDGM}.
In that case we can give a good description of the basin of attraction
of the quasi-periodic solution; cf. Section~\ref{sec:5}.
This provides a complementary result to the analysis performed
in~\cite{BDGM}, where the set of initial data generating
unbounded solutions was extensively studied.
The analysis can be extended to any function $x^{2p}$, $p\in\NNN$.
We also construct a positively invariant set containing
the attracting periodic orbit which allows us, together with the
results of Section~\ref{sec:5}, to obtain a larger set estimating
the basin of attraction (with area growing linearly in $\g$);
cf. Theorem~\ref{thm:5} in Section~\ref{sec:6}.

Finally we show that in the case $g(x)=x^{2p}$, $p\in\NNN$,
there are unbounded solutions which blow up in finite time, and we
discuss the implications for the varactor equation in~\cite{BDGM};
cf. Theorem~\ref{thm:6} in Section~\ref{sec:7}.

More formal statements of the results will be formulated
in the forthcoming sections. Some open problems will
be discussed at the end. Here we confine ourselves to
noting that, while in the case of periodic forcing
standard techniques, like those based
on Poincar\'{e} sections~\cite{L,M}, could be applied,
this is not the case for quasi-periodic forcing,
where no Poincar\'e maps can be introduced.

The rest of the paper is organised as follows.
Sections \ref{sec:2} and \ref{sec:3} are devoted to the
global study of the system~(\ref{eq:1.1}) for $\g$ large enough,
whereas in Section~\ref{sec:4} we draw the conclusions,
and we pass to the study of the system~(\ref{eq:1.4}).
In Section~\ref{sec:5} we specify the construction
envisaged for any $g(x)$ in~(\ref{eq:1.4}) to the case
$g(x)=x^{2p}$, with particular care for $p=1$, and
we look for an estimate of the basin of attraction of
the attracting periodic solution. Then, in Section~\ref{sec:6},
we use the techniques of~\cite{BDGM} to improve the estimate
of the basin of attraction. In Section~\ref{sec:7}
we show that, again in the case $g(x)=x^{2p}$,
there are unbounded solutions which blow up in finite time.
Finally in Section~\ref{sec:6} we mention some open problems,
and possible future directions of research.

\zerarcounters
\section{Existence of the quasi-periodic solution}
\label{sec:2}

First we show that for $\g$ large enough
there exists a quasi-periodic solution $x_{0}(t)$.

\begin{theorem}\label{thm:1}
Consider the equation~(\ref{eq:1.1}), with $f$ a non-zero average
quasi-periodic function analytic in its argument and with $\oo$
satisfying the Diophantine condition~(\ref{eq:1.3}).
There exists $\g_{0}>0$ such that for all $\g>\g_{0}$
there is a quasi-periodic solution $x_{0}(t)$ with the same
frequency vector as the forcing term. Such a solution extends to a
function analytic in $1/\g$ in a disc $\DD$ of the complex plane
tangent to the imaginary axis at the origin and centered on the
real axis. Furthermore, $x_{0}(t)=\al + O(1/\g)$,
where $\al=f_{\vzero}^{1/(2p+1)}\neq0$.

\end{theorem}

\prova We can apply the results of Section 7 in~\cite{GBD}.
If we set $g(x)=x^{2p+1}$,
then $g(c_{0})=f_{\vzero}$ yields $c_{0}=f_{\vzero}^{1/(2p+1)}$,
so that $g'(c_{0}) \neq 0$ as by assumption one has $f_{\vzero}\neq0$.
Both the existence of the analyticity domain $\DD$ and the form of
the solution itself follow from the analysis in~\cite{GBD}.$\EP$

If $\g$ is large enough, say $\g>\overline{\g}_{0} \ge \g_{0}$, then the
solution $x_{0}(t)$ is of definite sign. In the following we shall assume
that this is the case: hence $x_{0}(t) \neq 0$ for all $t\in\RRR$.

\zerarcounters
\section{Convergence to the quasi-periodic solution}
\label{sec:3}

Given the quasi-periodic solution $x_{0}(t)$ one can write
$x(t)=x_{0}(t)+\xi(t)$, with $\xi(t)$ satisfying the differential 
equation
\begin{equation}
\ddot \xi + \g \dot \xi + \xi F(\xi,x_{0}(t)) = 0  ,
\label{eq:3.1} \end{equation}
where we have defined
\begin{equation}
F(\xi,x) := \frac{1}{\xi} \left( (x+\xi)^{2p+1}-x^{2p+1} \right) =
\sum_{j=0}^{2p} \left( \begin{matrix} 2p+1 \\ j \end{matrix} \right)
\xi^{2p+1-j} x^{j} .
\label{eq:3.2} \end{equation}

We can write~(\ref{eq:3.1}) as
\begin{equation}
\begin{cases} \dot \x = y , \\
\dot y = - \g y - \xi F(\xi,x_{0}(t)) , \end{cases}
\label{eq:3.3} \end{equation}
that is $\dot z = \Phi(z)$, if we define $z=(\xi,y)$
and $\Phi(z) = (y,-\g y - \xi F(\xi,x_{0}(t))$.
We denote by $\f(t,z_0)$ the solution of~(\ref{eq:3.3})
with initial datum $z_0$.
Define also $P(\xi,t) := F(\xi,x_{0}(t))$ and $Q(\xi) := F(\xi,\al)$
and set $R(\xi,t) := P(\xi,t)/Q(\xi)$.

Here we prove the following result.

\begin{theorem}\label{thm:2}
Consider the equation~(\ref{eq:3.1}), with $x_{0}(t)$
the quasi-periodic solution of~(\ref{eq:1.1}) given in Theorem 1.
There exists $\g_{1}>0$ such that for all
$\g>\g_{1}$ all trajectories in phase space converge
toward the origin as time goes to infinity.
\end{theorem}

The proof will pass through several lemmata.

\begin{lemma}\label{lem:1}
Assume $\g>\overline{\g}_{0}$ so that $x_0(t)$ exists
and $x_{0}(t) \neq 0$ for all $t\in\RRR$.
There exist two positive constants $R_{1}$ and $R_{2}$ such that
\begin{equation}
R_{1} < R(\xi,t) < R_{2}
\label{eq:3.4} \end{equation}
for all $\xi\in\RRR$ and for all $t \in \RRR$.
\end{lemma}

\prova By~(\ref{eq:3.2}) we can write
\begin{equation}
F(\xi,x) = (2p+1) \int_{0}^{1} {\rm d} s \left(
x + s \, \xi \right)^{2p} ,
\label{eq:3.5} \end{equation}
so that $F(\xi,x) \geq 0$ for all $(\xi,x)\in\RRR^{2}$.
Moreover $F(0,0)=0$ and $F(\xi,x) >0$ for all $\xi\in\RRR$ if $x \neq 0$,
and $\lim_{|\xi|\to\infty}F(\xi,x)=\infty$ for all $x\in\RRR$.
Hence for $\al \neq 0$ and $\g>\overline{\g}_{0}$,
one has both $P(\xi,t)>0$ and $Q(\xi)>0$, hence also $R(\xi,t)>0$
for all $(\xi,t)\in\RRR^{2}$. Moreover
$\lim_{|\xi|\to\infty}R(\xi,t)=1$ for all $t\in\RRR$,
so that the assertion follows.$\EP$

\begin{lemma}\label{lem:2}
Consider the equation~(\ref{eq:3.1}), with $x_{0}(t)$
the quasi-periodic solution of~(\ref{eq:1.1}) given in Theorem 1.
There exists $\g_{2}>0$ such that for all $\g>\g_{2}$ there is
a convex set $\SSS$ containing the origin such that any trajectory
starting inside $\SSS$ tends to the origin as time goes to infinity.
One can take $\SSS$ such that $\partial\SSS$ crosses the positive and
negative $y$- and $\xi$-axes at distances $O(\g^{2})$
and $O(\g^{2/(p+1)})$ from the origin respectively.
\end{lemma}

\prova Rescale time through the Liouville transformation
\begin{equation}
\tau = \int_{0}^{t} {\rm d} t' \, \sqrt{R(\xi(t'),t')} ,
\label{eq:3.6} \end{equation}
which is well-defined by Lemma~\ref{lem:1}. Then,
if we introduce the coordinate transformation $\psi:(\xi,y)\to(v,y)$
by setting $\xi(t)=v(\tau(t))$ and
$y(t)=\sqrt{R(\xi(t),t)}\,w(\tau(t))$,
equation~(\ref{eq:3.3}) is transformed into
\begin{equation}
\begin{cases} v' = w , \\
w'= \displaystyle{
- \frac{w}{\sqrt{R}} \left( \g + \frac{R'}{2\sqrt{R}} \right) -
v \, Q(v) } , \end{cases}
\label{eq:3.7} \end{equation}
where primes denote differentiation with respect to $\tau$,
$Q(v(\tau))=Q(\xi(t(\tau))$
and $R=R(v(\tau),t(\tau))=R(\xi(t(\tau)),t(\tau))$.

The autonomous system
\begin{equation}
\begin{cases} v' = w , \\
w'= - v \, Q(v)  , \end{cases}
\label{eq:3.8} \end{equation}
can be explicitly solved: all trajectories move on the level curves
of the function
\begin{equation}
H(v,w) = \frac{1}{2} w^{2} + \int_{0}^{v} {\rm d} v' \, v' Q(v') .
\label{eq:3.9} \end{equation}
In~(\ref{eq:3.7}) one has $R'/\sqrt{R}=\dot R/R$, with
$$ \frac{\dot R}{R} = \frac{\dot P}{P} - \frac{\dot Q}{Q} , $$
and it is easy to see (Appendix \ref{app:A})
that there are two $\g$-independent
positive constants $B_{1}$ and $B_{2}$ such that
\begin{equation}
\left| \frac{\dot R}{2R} \right| <
\frac{1}{\g} \left( B_{1} + B_{2} |w| \right) .
\label{eq:3.10} \end{equation}
If $\g$ satisfies $\g^{2}>2B_{1}$ we can define $\widetilde w$ as
\begin{equation}
\widetilde w = \frac{\g^{2}-B_{1}}{B_{2}} > \frac{\g^{2}}{2B_{2}} , 
\label{eq:3.11} \end{equation}
so that $\overline \g := (\g+R'/2\sqrt{R})/\sqrt{R}>0$
for $|w| \leq \widetilde w$.

Consider the compact set $\widetilde\PPP$ whose
boundary $\partial\widetilde\PPP$ is given by the level curve
$H(v,w)=\widetilde w^{2}/2$ of the system~(\ref{eq:3.8}).
Such a curve crosses the $w$-axis at $w=\pm\widetilde w =
O(\g^{2})$ and the $v$-axis at $v=O(\g^{2/(p+1)})$.
If we take an initial datum $(v(0),w(0))\in\widetilde\PPP$
then the dissipation coefficient $\overline \g$ in~(\ref{eq:3.7}),
even if it changes with time, always remains strictly positive.
Moreover $H' = - \overline \g w^{2} \leq 0$ and $H'=0$ only for $w=0$,
and for $w=0$ the vector field in~(\ref{eq:3.7}) vanishes only at $v=0$,
because $Q(v) >0$ for all $v$ (cf. the proof of Lemma \ref{lem:1}).
Then we can apply Barbashin-Krasovsky theorem~\cite{BK,K},
and conclude that the origin is asymptotically stable
and that $\widetilde \PPP$ belongs to its basin of attraction.

Let $\PPP(t)$ be the time-dependent preimage of $\widetilde\PPP$ under
the transformation $\psi$. By Lemma~\ref{lem:1} if $\g$
is large enough there is a compact set $\SSS \subset \PPP(t)$
for all $t\in\RRR$, such that the boundary $\partial\SSS$ crosses
the positive and negative $y$- and $\xi$-axes at a distances of order $\g^{2}$
and $\g^{2/(p+1)}$ from the origin, respectively.
All trajectories starting from points inside $\SSS$
are attracted by the origin.$\EP$

\begin{lemma}\label{lem:3}
Consider the curve $g(\xi,t) = - \g^{-1} \xi F(\xi,x_{0}(t))$
in the plane $(\xi,y)$. There exists $\g_{3}>0$ such that for
$\g>\g_{3}$, outside the set $\SSS$
defined in Lemma~\ref{lem:2}, one has
\begin{equation}
- \frac{1}{2\g} \xi^{2p+1} \geq g(\xi,t) \geq - \frac{2}{\g} \xi^{2p+1}
\label{eq:3.12} \end{equation}
for all $t\in\RRR$.
\end{lemma}

\prova Consider $\xi\geq 0$ (the case $\xi<0$ can be discussed in the
same way). By~(\ref{eq:3.2}) one has
$$ \xi F(\xi,x_{0}(t)) = \xi^{2p+1} +
\sum_{j=1}^{2p} \left( \begin{matrix} 2p+1 \\ j \end{matrix} \right)
\xi^{2p+1-j} x_{0}^{j}(t) , $$
and, if $\g$ is sufficiently large so that $|x_{0}(t)| < 2|\al|$,
then for $\xi \geq 2|\al|$ the sum can be bounded
by $2^{2p+1}(2|\al|)\,\xi^{2p} \equiv C_{p}\xi^{2p}$. Hence one has
\begin{equation}
\frac{1}{2}\,\xi^{2p+1} \leq \xi F(\xi,x_{0}(t))
\leq 2 \xi^{2p+1} , 
\label{eq:3.13} \end{equation}
as soon as $\xi \geq 2C_{p}$ (note that if $\xi \geq 2C_{p}$ then one
has automatically $\xi \geq 2|\al|$). Next, we want to show
that the latter inequality is satisfied outside $\SSS$.

\begin{figure}[htbp]
\begin{centering}
\ins{160pt}{-135pt}{$y_{0}$}
\ins{246pt}{-22pt}{$\xi_{0}$}
\ins{164pt}{-6pt}{$y$}
\ins{304pt}{-40pt}{$\xi$}
\ins{280pt}{-4pt}{$y_{1}(\xi)$}
\ins{248pt}{-160pt}{$y_{2}(\xi)$}
\ins{213pt}{-20pt}{$\overline\xi$}
\ins{256pt}{-60pt}{$\partial\SSS$}
\includegraphics*[angle=0,width=2.4in]{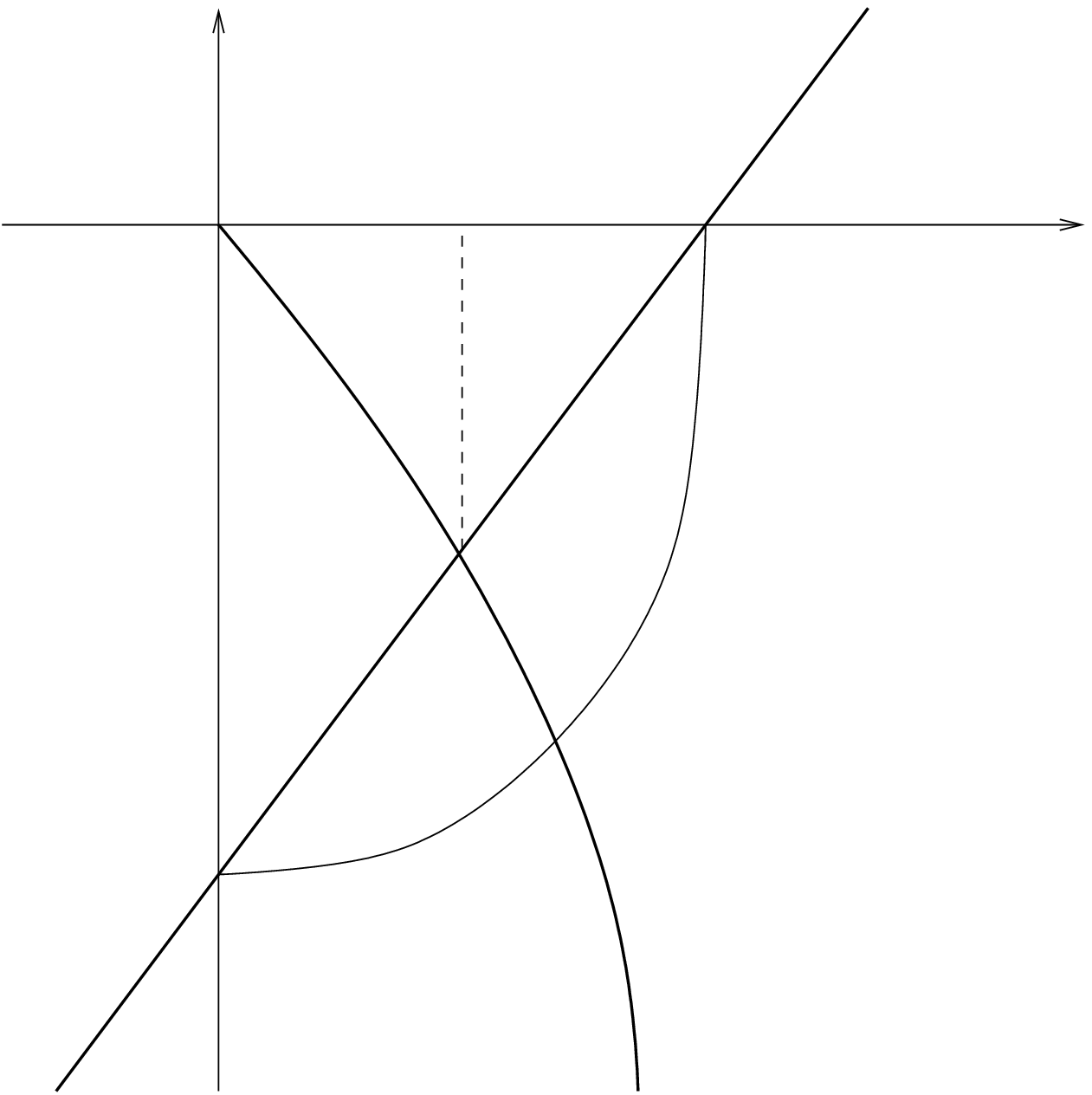}
\caption{\label{fig3-1}
Construction used in the proof of Lemma \ref{lem:3}.}
\end{centering}
\end{figure}

Consider the intersection of the graph of $g(\xi,t)$ with $\partial\SSS$.
Let II be the quadrant $\{(\xi,y) \in \RRR^{2} : \xi \geq 0 , y < 0 \}$;
cf. figure~\ref{fig3-1}. The curve $\partial\SSS \cap$II is below the line
\begin{equation}
y_{1}(\xi) = y_{0} \left( 1 - \frac{\xi}{\xi_{0}} \right) ,
\label{eq:3.14} \end{equation}
where $y_{0}:=-b \g^{2}$, with $b>0$, is the $y$-coordinate of the point
at which $\partial\SSS$ crosses the $y$-axis, and
$\xi_{0}:=a \g^{2/(p+1)}$, with $a>0$, is the $\xi$-coordinate of
the point at which $\partial\SSS$ crosses the $\xi$-axis.
On the other hand the graph of $g(\xi,t)$ in II is above the curve
\begin{equation}
y_{2}(\xi) = -\frac{2}{\g} \left( \xi + 2|\al| \right)^{2p+1} ,
\label{eq:3.15} \end{equation}
because in (\ref{eq:3.2}) one has $\xi F(\xi,x) \le
|x+\xi|^{2p+1} + |x|^{2p+1}$.

As a consequence in II the
two curves $\partial\SSS$ and $g(\xi,t)$ cannot 
cross each other for $\xi\in(0, a\g^{\be}/2]$, with $\be\leq 2/(2p+1)$.
The latter assertion can be proved by \textit{reductio ad absurdum}.
First note that $a\g^{\be} \leq \xi_{0}$ for $\g$ large enough.
Suppose that there exists $\overline\xi \leq a \g^{\be}/2$ such that
$y_{1}(\overline\xi)=y_{2}(\overline\xi)$. Then one has
$$ b \frac{\g^{2}}{2} \le b \g^{2} \left( 1 - 
\frac{\overline\xi}{\xi_{0}} \right)
= \frac{2}{\g} \left( \overline\xi + 2|\al| \right)^{2p+1}
\le \frac{4}{\g} \max\{\overline\xi, 2|\al|\}^{2p+1} , $$
that is $b \g^{3} \leq 8 \max\{\overline\xi, 2|\al|\}^{2p+1}$,
which is not possible if $\be\leq 2/(2p+1)$ and $\g$ is large enough.

Therefore in II the graph of $g(\xi,t)$ can be outside $\SSS$ only
for $\xi > a\g^{\be}/2$, which is greater than $2C_{p}$
for $\g$ large enough. Hence~(\ref{eq:3.13}) is satisfied outside 
$\SSS$, so that~(\ref{eq:3.12}) follows. $\EP$

\begin{lemma}\label{lem:4}
Consider the equation~(\ref{eq:3.1}), with $x_{0}(t)$
the quasi-periodic solution of~(\ref{eq:1.1}) given in Theorem 1.
There exists $\g_{4}>0$ such that for all $\g>\g_{4}$,
if $z \notin \SSS$,then either $\f(t,z)$ enters $\SSS$
or crosses the $y$-axis outside $\SSS$ in a finite positive time.
\end{lemma}

\prova Consider the four quadrants
\begin{eqnarray}
{\rm I} & = & \{ (\xi,y) \in \RRR^{2} :
\xi > 0 , y \geq 0 \} , \nonumber \\
{\rm II} & = & \{ (\xi,y) \in \RRR^{2} :
\xi \geq 0 , y < 0 \} , \nonumber \\
{\rm III} & = & \{ (\xi,y) \in \RRR^{2} :
\xi < 0 , y \leq 0 \} , \nonumber \\
{\rm IV} & = & \{ (\xi,y) \in \RRR^{2} :
\xi \leq 0 , y > 0 \} .
\label{eq:3.16} \end{eqnarray}
In I one has $\dot \xi \geq 0, \dot y <0$,
in II one has $\dot \xi < 0$,
in III one has $\dot \xi \leq 0, \dot y >0$,
and in IV one has $\dot \xi > 0$.
It is easy to see that each trajectory starting in I
enters II and each trajectory starting from III enters IV
in a finite time (see Appendix \ref{app:B}).

\begin{figure}[htbp]
\begin{centering}
\ins{165pt}{-77pt}{$\overline y$}
\ins{245.5pt}{-20pt}{$\overline \xi$}
\ins{165pt}{-6pt}{$y$}
\ins{304pt}{-40pt}{$\xi$}
\ins{282pt}{-100pt}{$\CC_{1}$}
\ins{256pt}{-60pt}{$\partial\SSS$}
\ins{200pt}{-100pt}{$\SSS$}
\ins{290pt}{-70pt}{$\TT_{1}$}
\ins{220pt}{-150pt}{$\TT_{2}$}
\includegraphics*[angle=0,width=2.4in]{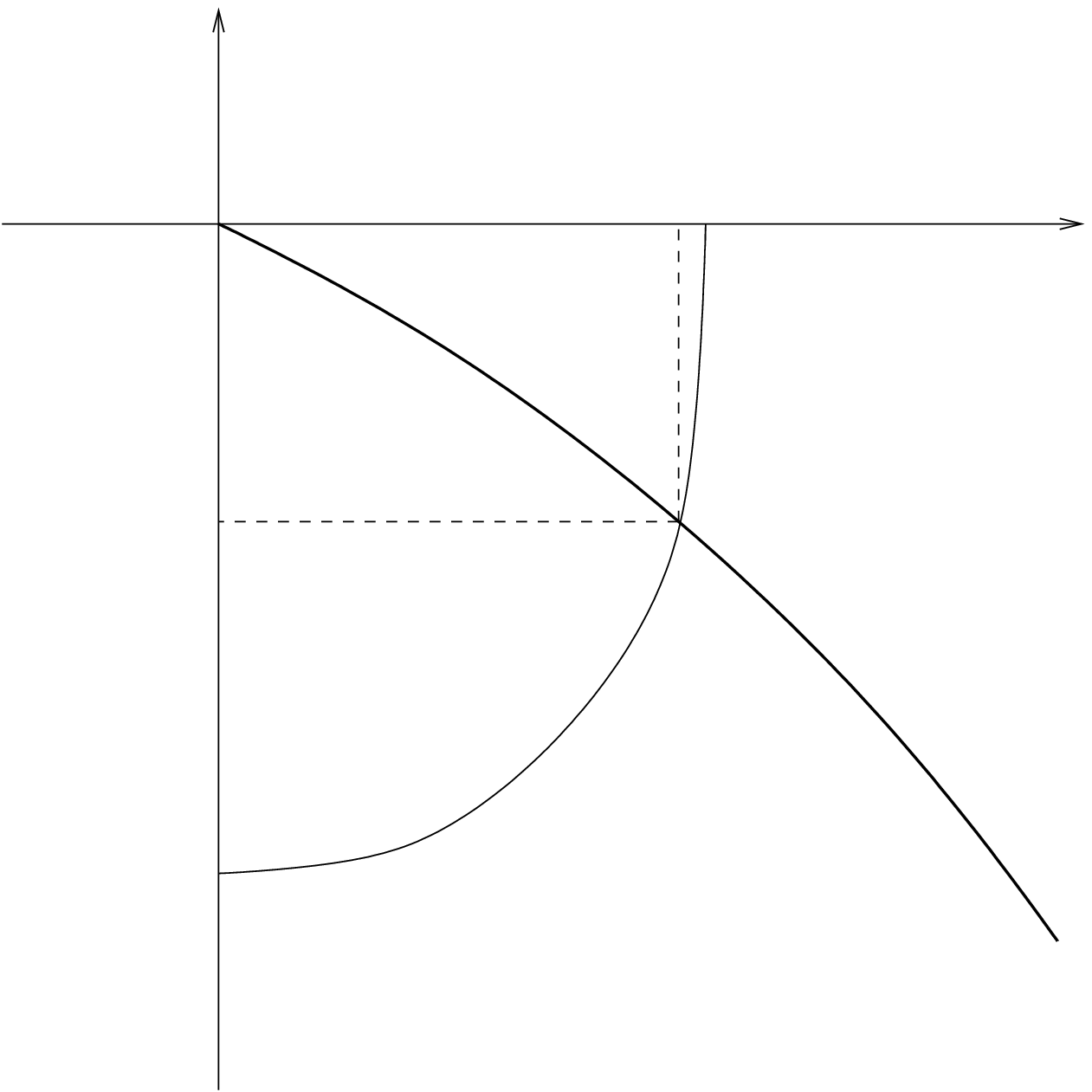}
\caption{\label{fig3-2}
Construction used in the proof of Lemma \ref{lem:4}.}
\end{centering}
\end{figure}

Consider now an initial datum $z$ in II but not in to $\SSS$.
Let $\CC_{1}$ be a continuous curve $\xi \to y(\xi)$ in II such that
$\dot y < 0$ for $z$ in II above $\CC_{1}$; cf. figure~\ref{fig3-2}.
Existence of such a curve follows from Lemma~\ref{lem:3},
which also implies that $\CC_{1}$ is decreasing outside $\SSS$
(see Appendix \ref{app:C}). The curve $\CC_{1}$ divides II
into two sets IIa and IIb, with IIa above IIb.
Denote by $\TT_{1}$ and $\TT_{2}$ the parts of IIa and IIb,
respectively, outside $\SSS$. Hence for $z \in \TT_{1}$
the trajectory $\f(t,z)$ either enters $\SSS$ or enters $\TT_{2}$.
In the latter case it cannot come back to $\TT_{1}$,
hence $y(t) \leq \overline{y}$, if $(\overline \xi,\overline y)=
\CC_{1} \cap \partial\SSS$. This means that if the solution does not enter
$\SSS$ then it has to cross the vertical axis and enter III.

Analogously one discusses the case of initial data $z$ in IV,
outside $\SSS$: their evolution leads either to $\SSS$
or to I. Hence the lemma is proved.$\EP$

\begin{lemma}\label{lem:5}
Consider the equation~(\ref{eq:3.1}), with $x_{0}(t)$ the
quasi-periodic solution of~(\ref{eq:1.1}) given in Theorem 1.
There exists $\g_{5}>0$ such that for all $\g>\g_{5}$ and
for all $z\notin\SSS$ on the vertical axis, either $\f(t,z)$
enters $\SSS$ or the trajectory $\f(t,z)$ re-crosses the vertical axis
at a point $z_{1}$ which is such that $|z|-|z_{1}| > \delta$
for some positive $z$-independent constant $\delta$.
\end{lemma}

\begin{figure}[htbp]
\begin{centering}
\ins{209pt}{-6pt}{$Y$}
\ins{326pt}{-150pt}{$X$}
\ins{210pt}{-28pt}{$1$}
\ins{248pt}{-212pt}{$\CC_{2}$}
\ins{266pt}{-106pt}{$\partial S$}
\ins{234pt}{-122pt}{$S$}
\ins{202pt}{-262pt}{$-1$}
\ins{260pt}{-232pt}{$P$}
\ins{280pt}{-250pt}{$Q$}
\ins{190pt}{-222pt}{$Y(T_{3})$}
\ins{310pt}{-60pt}{$\Gamma$}
\ins{292pt}{-31pt}{$(X(t),Y(t))$}
\includegraphics*[angle=0,width=3.2in]{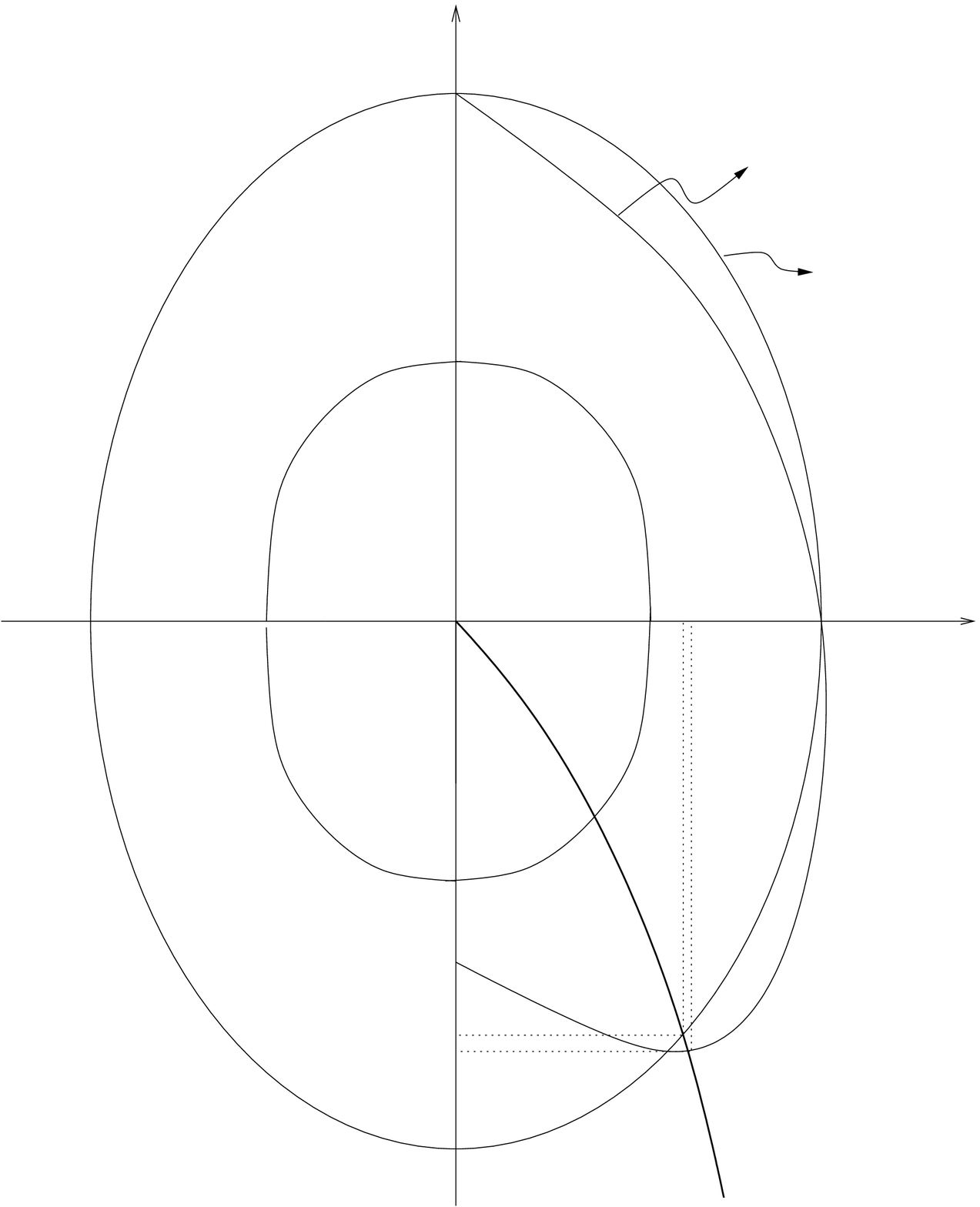}
\caption{\label{fig3-3}
Construction used in the proof of Lemma \ref{lem:5}.}
\end{centering}
\end{figure}

\prova Fix an initial datum $z\notin\SSS$ on the vertical axis.
This means that at $t=0$ one has $z=z(0)=(0,y(0))$ outside $\SSS$.
Assume for concreteness $y(0)>0$ and set $y(0)=1/\eps^{p+1}$,
with $\eps>0$. As $z(0)\notin\SSS$ there exists a constant $C_{2}$ such
that $\g^{2}\eps^{p+1} \le C_{2}$. Consider the change of coordinates
\begin{equation}
X = \eps \xi , \qquad Y = \eps^{p+1} y , \qquad
T = \frac{t}{\eps^{p}} .
\label{eq:3.17} \end{equation}
Then the system~(\ref{eq:3.3}) becomes
\begin{equation}
\begin{cases} X' = Y , \\
Y' = - \g \eps^{p} Y - X F(X,\eps x_{0}(\eps^{p} T) , \end{cases}
\label{eq:3.18} \end{equation}
where primes denote differentiation with respect to $T$.
Note that $Y(0)=1$ and $X F(X,\eps x_{0}(\eps^{p}T)=X^{2p+1}+
O(\eps \, X^{2p})$. Call $S$ the image of $\SSS$
under the transformation (\ref{eq:3.17}); cf. figure~\ref{fig3-3}.

We can rewrite the system~(\ref{eq:3.18}) as
\begin{equation}
\begin{cases} X' = Y , \\
Y' = \Psi(X,Y) \equiv
\Psi_{1}(X,Y) + \Psi_{2}(X,Y) + \Psi_{3}(X,Y) , \end{cases}
\label{eq:3.19} \end{equation}
with
\begin{equation}
\begin{cases}
\Psi_{1}(X,Y)=-X F(X, \eps x_{0}(0)) , \\
\Psi_{2}(X,Y)=- X \left[\;
F(X, \eps x_{0}(\eps^{p} T)) - F(X,\eps x_{0}(0)) \;\right] , \\
\Psi_{3}(X,Y)=-\g \eps^{p} Y . \end{cases}
\label{eq:3.20} \end{equation}

If we replace $\Psi(X,Y)$ with $\Psi_{1}(X,Y)$ in~(\ref{eq:3.19})
the trajectory moves on the level curve
$\Gamma=\{(X,Y)\in\RRR^{2} : H(X,Y)=1/2\}$ for the function
\begin{equation}
H(X,Y) = \frac{1}{2} Y^{2} + \int_{0}^{X}
{\rm d} X \, X F(X, \eps x_{0}(0)) 
\label{eq:3.21} \end{equation}
and crosses the vertical axis at the point $(0,-1)$, hence 
at the same distance from the origin as at $t=0$.
By Lemma~\ref{lem:3} there exists in II a curve $\CC_{2}$, decreasing
outside $S$, such that
$$ -\g\eps^{p} Y \geq \Psi(X,Y) \ge -\g\eps^{p} Y/2 , $$
for $Y$ in II below $\CC_{2}$ (see Appendix \ref{app:C}).
Such a curve can be chosen in such a way that it crosses
the level curve $\Gamma$ in a point $P=(X_{P},Y_{P})$,
with $X_{P}=2D_{1} (\g\eps^{p})^{1/(2p+1)}$, for some constant 
$D_{1}$ (see Appendix \ref{app:C}). Note that the time $T_{1}$
necessary to reach such a point is of order 1.

If we take into account the component $\Psi_{2}(X,Y)$
of the vector field in~(\ref{eq:3.19}), we can move from
$P$ at most by a quantity of order $\eps^{p+1}$. Indeed, as long as
the motion remains close to that generated by the vector field
$\Psi_{1}(X,Y)$, one has
\begin{equation}
\left| F(X, \eps x_{0}(\eps^{p} T_{1})) - F(X,\eps x_{0}(0)) \right|
\leq D_{2}' \eps \left| x_{0}(\eps^{p} T_{1}) - x_{0}(0) \right|
\leq D_{2}'' \eps^{p+1} ,
\label{eq:3.22} \end{equation}
for suitable positive constants $D_{2}'$ and $D_{2}''$,
so that the points reached at time $T_{1}$ by moving
according to the vector fields $\Psi_{1}$ and $\Psi_{1}+\Psi_{2}$
cannot be more distant than $D_{2}\eps^{p+1}$ for some constant $D_{2}$.
This follows from the fact that the system is quasi-integrable,
so that in a time of order 1 the action variable can change at most
by a quantity of order of the perturbation as bounded in~(\ref{eq:3.22});
see Appendix \ref{app:D}.

Finally the component $\Psi_{3}(X,Y)$ points inward along the
full length of the curve $\Gamma$.
Define $T_{2}$ as the time at which the
the trajectory of the full system~(\ref{eq:3.19}) crosses
the curve $\CC_{2}$ in a point $Q$ near $P$. Of course
$T_{2}$ is near $T_{1}$, and so is of order 1, and $X(T_{2}) \geq X_{P}/2$
by construction, while $Y(T_{2}) \geq Y_{P} - D_{2}\eps^{p+1}$,
with $Y_{P} > - 1$.

From time $T_{2}$ onwards, we have
\begin{equation}
\begin{cases} X' = Y , \\
Y' \geq -\g\eps^{p} Y / 2 , \end{cases}
\label{eq:3.23} \end{equation}
as long as the motion remains below $\CC_{2}$. The latter property
is easily checked to hold (see Appendix \ref{app:C}).
Then the trajectory crosses the vertical axis and
meanwhile, at least, moves upward in the vertical direction by a
quantity $\g\eps^{p} X_{P}/2 = D_{1} \g\eps^{p}(\g\eps^{p})^{1/(2p+1)}$.

Therefore when the trajectory again crosses the vertical axis,
this happens at a time $T_{3}$ such that
$Y(T_{3}) \geq Y_{P} - D_{2} \eps^{p+1} + D_{1}
\g\eps^{p}(\g\eps^{p})^{1/(2p+1)} > - 1 + \Delta Y$, with
$$ \Delta Y = D_{1} \g\eps^{p}(\g\eps^{p})^{1/(2p+1)} - D_{2} \eps^{p+1}
\geq D_{2} \eps^{p+1} , $$
where the latter inequality holds provided $(\g\eps^{p})^{1+1/(2p+1)}
\geq 2D_{2} \eps^{p+1}$, that is provided
\begin{equation}
\g^{2(p+1)} \geq D_{0} \eps^{p+1} , \qquad 
D_{0} = (2D_{2})^{2p+1} .
\label{eq:3.24} \end{equation}
Since $\eps^{p+1}\g^{2} \leq C_{2}$, inequality~(\ref{eq:3.24})
is satisfied if $\g^{2(p+1)} \geq D_{0} C_{2} \g^{-2}$,
which requires $\g^{2(p+2)} \geq D_{0} C_{2}$, that is
\begin{equation}
\g \geq \left( D_{0} C_{2} \right)^{1/2(p+2)} , \qquad
D_{0} = (2D_{2})^{2p+1} .
\label{eq:3.25} \end{equation}
Under this condition one has $|Y(0)|-|Y(T_{3})| = 1 - |Y(T_{3})| \geq
\Delta Y \geq D_{2} \eps^{p+1}$, so that,
in terms of the original coordinate $y$, one has
$|y(0)|-|y(t_{3})| \geq D_{2}$.

Then, if the trajectory crosses the
vertical axis once more in the positive direction
(and this necessarily happens if it does not enter $\SSS$,
by Lemma~\ref{lem:4}), this occurs at a time $t_{4}$ such that
\begin{equation}
|y(0)|-|y(t_{4})| \geq 2 D_{2} ,
\label{eq:3.26} \end{equation}
where we recall that the constant $D_{2}$ is independent of the
initial datum $y(0)$. Simply one can repeat the argument above
by taking $(0,y(t_{3}))$ as initial datum and calling $t_{4}$
the time of crossing of the positive $\xi$-axis. This means that the
trajectory either enters $\SSS$ or, after a complete cycle, moves closer
to the origin by a finite positive distance $\delta=2D_{0}$.$\EP$

\begin{lemma}\label{lem:6}
Consider the equation~(\ref{eq:3.1}), with $x_{0}(t)$
the quasi-periodic solution of~(\ref{eq:1.1}) given as in Theorem 1.
There exists $\g_{6}>0$ such that for all $\g>\g_{6}$ for all $z\notin\SSS$
there is a finite time $t(z)$ such that $\f(t(z),z) \in \SSS$. 
\end{lemma}

\prova Consider $z\notin\SSS$. By Lemma~\ref{lem:4} either
$\f(t,z)$ enters $\SSS$ or there exists a time $t_{1}$ such that
$\f(t_{1},z)$ is on the vertical axis outside $\SSS$.
Hence, without loss of generality, we can consider only
initial data $z=(\xi,y)$ outside $\SSS$ such that $\xi=0$.
Assume $y>0$ (if $y<0$ the discussion proceeds in the same way):
we can apply Lemma~\ref{lem:4} and we find that, 
as far as the trajectory does not enter $\SSS$,
at each turn it gets closer to $\SSS$ by a finite quantity.
Hence sooner or later it enters $\SSS$.$\EP$

Theorem~\ref{thm:2} follows from the lemmata above: it is enough
to take $\g_{1}=\max\{\overline{\g}_{0},\g_{2},\g_{3},\g_{4},
\g_{5},\g_{6}\}$, so that all lemmata apply.

\zerarcounters
\section{Uniqueness of the quasi-periodic solutions}
\label{sec:4}

Let $\g$ be $\g>\max\{\g_{0},\g_{1}\}$. Then there exists
a quasi-periodic solution $x_{0}(t)$ for the system~(\ref{eq:1.1}),
and such a solution is a global attractor.

In~\cite{GBD} we explicitly constructed a solution $x_{0}(t)$
with the properties stated in Theorem~\ref{thm:1}.
Such a solution turns out to be Borel summable for $d=1$.
In general the solution is obtained from the formal series 
through a suitable summation procedure. Since Theorem~\ref{thm:2}
implies that if there exists a quasi-periodic solution $x_{0}(t)$
this has to be unique, we can conclude that for real $\g$ large enough
there exists a unique quasi-periodic solution $x_{0}(t)$ with the
same frequency vector as the forcing $f$. By setting $\eps=1/\g$,
in the complex $\eps$-plane, there is a solution $x_{1}(t)$
which is analytic in a domain containing
a disk $\DD$ with centre on the real axis
and tangent at the origin to the imaginary axis.
For real $\g$ such a solution coincides with $x_{0}(t)$
(as the latter is the only one), hence, by uniqueness of the
analytic continuation, the function $x_{1}(t)$ is the only
solution of~(\ref{eq:1.1}) in all the domain $\DD$.
In particular it is the only one which admits the
formal expansion given by perturbation theory.

We can summarise the discussion above through the following statement.

\begin{theorem}\label{thm:3}
Consider the equation~(\ref{eq:1.1}), with $f$ a non-zero average
quasi-periodic function analytic in its argument and
with $\oo$ satisfying the Diophantine condition~(\ref{eq:1.3}).
There exists $\g_{0}>0$ such that for all real
$\g>\g_{0}$ there is a unique quasi-periodic solution $x_{0}(t)$ with the
same frequency vector as the forcing term. Such a solution describes a
limit cycle in the plane $(x,\dot x)$ which is a global attractor.
\end{theorem}

Note that the hypotheses made in Theorem~\ref{thm:1} are more
restrictive than those considered in~\cite{GBD}. In particular
we excluded both polynomial nonlinearities and monomial
nonlinearities with even degree.

The first request aims to have a strictly positive
function $R(\xi,t)$, which was used in order to construct the
positively invariant set $\SSS$. We leave as an open problem
to study what happens if the nonlinearity $x^{2p+1}$
in~(\ref{eq:1.1}) is replaced with
\begin{equation}
g(x) = \sum_{j=1}^{2p+1} a_{j} \, x^{j} , \qquad a_{j} \in \RRR ,
\qquad a_{2p+1} > 0 .
\label{eq:4.1} \end{equation}
For $d=1$ and $\g=0$ it is known that all motions are bounded, 
also replacing the constants $a_{j}$ with periodic
functions~\cite{DZ,LL,CNY}. The same holds for $d>1$~\cite{LZ}.
One could expect that the presence of friction tends to contract
phase space toward some periodic solution (which certainly exists
for $\g$ large enough, as proved in \cite{GBD}), but our results
do not allow us to treat, in general, such a case.

If the nonlinearity $x^{2p+1}$ is replaced
with an even monomial $x^{2p}$, with $p\in\NNN$,
then, under the further condition that $f_{\vzero}>0$,
there is a quasi-periodic solution $x_{0}(t)$:
again this follows from~\cite{GBD}. In such a case
$x_{0}(t)$ is not a global attractor, as there are unbounded
solutions; cf. for example~\cite{BDGM} for $p=1$.
Still one can prove that the solution found in~\cite{GBD}
is unique, in the sense that it is the only attractor
in a neighbour of the solution itself,
and it is a local attractor. The same result holds,
more generally, for any analytic $g(x)$ in~(\ref{eq:1.4})
such that $g(c_{0})=f_{\vzero}$ and $g'(c_{0})\neq0$
for some $c_{0}\in\RRR$. A more formal statement is as follows.

\begin{theorem}\label{thm:4}
Consider the equation~(\ref{eq:1.4}), with $f$,
given by~(\ref{eq:1.2}), and $g$ both analytic in their arguments,
and with $\oo$ satisfying the Diophantine condition~(\ref{eq:1.3}).
Assume that there exists $c_{0}\in\RRR$ such that $g(c_{0})=f_{\vzero}$
and $g'(c_{0})\neq0$. There exists $\g_{0}>0$ such that for $\g>\g_{0}$
there is a unique quasi-periodic solution $x_{0}(t)$
which has the same frequency vector as $f$, reduces to $c_{0}$
in the limit $\g\to\infty$, and extends to a function
analytic in a disk with center on the positive real axis
and boundary tangent to the vertical axis at the origin.
Furthermore, there exists $\g_{1} \geq \g_{0}$ such that
for $\g>\g_{1}$ there is a neighbourhood $\BB$ of the
point $(c_{0},0)$, containing the orbit described by $x_{0}(t)$,
with the property that all trajectories starting in $\BB$
are attracted to the cycle described by $x_{0}(t)$ in the plane.
\end{theorem}

\prova The existence of a quasi-periodic solution $x_{0}(t)$ with
the same frequency vector $\oo$ as the forcing was proved in~\cite{GBD}.
As a byproduct of the proof, one can write $x_{0}(t)=c_{0}+x_{1}(\oo t)$,
with $x_{1}(\pps)$ analytic in $\pps$ and of order $\eps$,
if $\eps=1/\g$ (that is for $\eps$ small enough one has
$|x_{1}(\pps)| \le C |\eps|$ for all $\pps$ and for a suitable $C$).
Therefore we can write $x(t)=x_{0}(t)+\xi(t)$, where $\xi(t)$
satisfies the differential equation
\begin{equation}
\ddot \xi + \g \dot \x + \xi F(\xi,x_{0}(t)) = 0 ,
\label{eq:4.2} \end{equation}
with
\begin{equation}
F(\xi,x)= \frac{1}{\xi} \Big( g(x+\xi) - g(x) \Big) =
\partial_{x} g (x) + O(\xi) .
\label{eq:4.3} \end{equation}
Then we can write~(\ref{eq:4.2}) as a system of
first order differential equations,
\begin{equation}
\begin{cases} \dot \x = y , \\
\dot y = - \g y - \xi F(\xi,x_{0}(t)) , \end{cases}
\label{eq:4.4} \end{equation}
and define $R(\xi,t)=F(\xi,x_{0}(t))/F(\xi,c_{0})$.
It is easy to see that one has
$$ \lim_{\xi\to 0} R(\xi,t) = 1 + O(\eps) , $$
so that for $\eps$ small enough one has $R_{1}<R(\xi,t)<R_{2}$,
for two suitable positive constants $R_{1}$ and $R_{2}$.

Then we can rescale time and variables by setting
\begin{equation}
\tau = \int_{0}^{t} {\rm d}t \, \sqrt{R(\xi(t'),t')} ,
\qquad \xi(t) = v(\tau(t)) , \qquad
\qquad y(t) = \sqrt{R(\xi(t),t)} \, w(\tau(t)) ,
\label{eq:4.5} \end{equation}
which transforms the system~(\ref{eq:4.4}) into
\begin{equation}
\begin{cases} v' = w , \\
w'= \displaystyle{
- \frac{w}{\sqrt{R}} \left( \g + \frac{R'}{2\sqrt{R}} \right) -
v \, F(v,c_{0}) } , \end{cases}
\label{eq:4.6} \end{equation}
where primes denote differentiation with respect to $\tau$.

If we neglect the friction term in~(\ref{eq:4.6}) we obtain
the autonomous system
\begin{equation}
\begin{cases} v' = w , \\
w'= - v F(v,c_{0})) = - \partial_{x} g(c_{0})\, v +
O(v^{2}) , \end{cases}
\label{eq:4.7} \end{equation}
which admits the constant of motion
\begin{equation}
H(v,w) = \frac{1}{2} w^{2} +
\frac{1}{2} \partial_{x} g(c_{0}) \, v^{2} + O(v^{3}) , \qquad
\partial_{x} g(c_{0}) \neq 0 .
\label{eq:4.8} \end{equation}
Hence the origin is a stable equilibrium point for~(\ref{eq:4.7}),
and the level curves for $H$ are close to ellipses in a neighbourhood
$\widetilde\PPP$ of the origin. It is easy to check that in
$\widetilde\PPP$ the coefficient of the friction term is strictly
positive, for $\g$ large enough, because $R'/2\sqrt{R}$ is
(in $\widetilde\PPP$) less than a constant.
Hence we can apply Barbashin-Krasovsky's theorem and conclude that
the origin is asymptotically stable and $\widetilde\PPP$
belongs to its basin of attraction. If we go back to the original
variables $(\xi,y)$ we find that $\widetilde\PPP$ is transformed back
to a time-dependent set $\PPP(t)$. But the dependence on $t$ of $\PPP(t)$
is very weak (as $R$ is close to 1 for $\g$ large enough), so that there
exists a convex set $\SSS \subset \PPP(t)$ for all $t\in\RRR$.
Hence any trajectory starting from $\SSS$ is attracted toward the origin.
In terms of the variables $(x,y)$, using once more
that the solution $x_{0}(t)$ is close to $c_{0}$ within $O(\eps)$,
we can say that, for $\eps$ small enough (that is for $\g$ large
enough) there exists a neighbourhood $\BB$ of the point $(c_{0},0)$
such that it contains the cycle described by the quasi-periodic
solution $x_{0}(t)$ in the plane $(x,\dot x)$, and
any trajectory starting from $\BB$ is attracted by such a cycle.

In particular the solution $x_{0}(t)$ is the only quasi-periodic
solution which tends to $c_{0}$ as $\eps\to 0$, and
for $\eps>0$ small enough, say $\eps<\eps_{0}$, it is the only one
which admits the formal power expansion~(\ref{eq:1.6}).
Such a solution was proved in~\cite{GBD} to be analytic in
a domain $\DD$ containing the interval $(0,\eps_{0})$,
hence by the uniqueness of the analytic continuation,
we can conclude that $x_{0}(t)$ is unique in all $\DD$.$\EP$

By looking at the proof of Theorem~\ref{thm:4}, we see that
it proceeds along the same lines of the proof of Lemma~\ref{lem:2}
in Section~\ref{sec:3}. By using the definitions of~\cite{ER},
Theorem~\ref{thm:4} says that the orbit described by $x_{0}(t)$
is an attracting set (and an attractor)
with fundamental neighbourhood $\BB$.

\zerarcounters
\section{The varactor equation, and generalisations}
\label{sec:5}

The construction of the set $\BB$ in the proof of Theorem~\ref{thm:4}
can be improved of course in concrete examples.
Here we want to consider explicitly the case of even monomials
$x^{2p}$, with particular emphasis for $p=1$.
Hence we discuss explicitly the case $p=1$,
but the discussion can be easily extended to all $p \ge 1$.

So, let us fix $p=1$. If we define the functions $F(\xi,x)$
and $R(\xi,t)$ as in the proof of Theorem~\ref{thm:4}, we find
\begin{equation}
\xi \, F(\xi,x) = \xi^{2} + 2 x \xi , \qquad
R(\xi,t) = \frac{\xi + 2x_{0}(t)}{\xi + 2 c_{0}} =
1 + 2 \frac{x_{0}(t) - c_{0}}{\xi + 2 c_{0}} ,
\label{eq:5.1} \end{equation}
where $c_{0}=\sqrt{f_{\vzero}}$, with $f_{\vzero}>0$.
For $\g$ large enough, that is for $\eps$ small enough, one has
\begin{equation}
|x_{0}(t) - c_{0}| < C_{1} \eps ,
\label{eq:5.2} \end{equation}
for a suitable constant $C_{1}$
(see the first lines of the proof of Theorem~\ref{thm:4}).

The transformation of coordinates described in the proof of
Theorem~\ref{thm:4} leads to~(\ref{eq:4.6}), with
\begin{equation}
v \, F(v,c_{0}) = v^{2} + 2 c_{0} v ,
\label{eq:5.3} \end{equation}
so that, if we set $\overline \g=(\g+R'/2\sqrt{R})/\sqrt{R}$,
we can interpret the system as a Hamiltonian system with hamiltonian
\begin{equation}
H(v,w) = \frac{1}{2} w^{2} + U(v) , \qquad
U(v) = \frac{1}{3} v^{3} + c_{0} v^{2} ,
\label{eq:5.4} \end{equation}
in the presence of a friction term with non-constant
dissipation coefficient $\overline\g$. If we neglect the friction term,
that is if we put $\overline\g=0$, then the system becomes
the Hamiltonian system
\begin{equation}
\begin{cases} v' = w , \\
w'= - v^{2} - 2 c_{0} v , \end{cases}
\label{eq:5.5} \end{equation}
which admits a stable equilibrium point $P_{1}=(0,0)$,
an unstable equilibrium point $P_{2}=(-2c_{0},0)$
and the separatrix $\Gamma$ of equation
\begin{equation}
w = \pm \sqrt{ 2 ( H(-2c_{0},0) - U(v)) } , \qquad
H(-2c_{0},0) = \frac{4c_{0}^{3}}{3} ,
\label{eq:5.6} \end{equation}
which contains a homiclinic orbit to the right of $P_{2}$.
Consider the open set $\DD_{1}$ containing the point $P_{1}$ and
with boundary given by (the closure of) the homoclinic orbit in $\Gamma$.
Any level curve with energy $V(0)<E<V(-2c_{0})$ contains
a bounded connected component, internal to $\DD_{1}$,
which is a closed orbit $\CC_{E}$ for the system~(\ref{eq:5.5}).
Call $\CC$ the closed curve $\CC_{E}$ such that
it intersects the negative $v$-axis at $v=-2c_{0}+C_{2} \eps^{\be}$,
for suitable constants $C_{2}$ and $\be$ to be fixed,
and $\DD$ the bounded open set with boundary $\CC$.

From~(\ref{eq:5.1}), one has 
\begin{equation}
\left| R(\xi,t) - 1 \right| < \frac{C_{1}}{C_{2}} \eps^{1-\be} ,
\label{eq:5.7} \end{equation}
so that $R(\xi,t) \to 1$ as $\eps \to 0$, provided $\be<1$;
for instance one can take $\be=1/2$. Then, if $\g$ is large enough,
for any point $(v,w)\in \overline{\DD}$ one has
\begin{equation}
\left| R'(\xi,t) \right| \le
\left| \frac{2 \xi' (c_{0}-x_{0})}{(\xi+2c_{0})^{2}} \right| +
\left| \frac{2 x_{0}'}{\xi+2c_{0}} \right| \le
C_{3} \eps^{1-2\be} + C_{4} \eps^{1-\be} ,
\label{eq:5.8} \end{equation}
for suitable constants $C_{3}$ and $C_{4}$, so that,
by using also that $\xi=v$ and
\begin{equation}
\left| \xi + 2c_{0} \right| \ge C_{2}\eps^{\be}
\label{eq:5.9} \end{equation}
for all $(v,w) \in \overline{\DD}$, one obtains
\begin{equation}
\left| \g-\g' \right| \le \left| \frac{R'}{2R} \right| \le
C_{5} \eps^{1-2\be} \le C_{5} \eps^{2(1-\be)} \g < \g
\label{eq:5.10} \end{equation}
for a suitable constant $C_{5}$, hence $\g'>0$ for $\g$ large enough.
Therefore we can conclude that $\overline{\DD}$ is positively invariant,
and we can apply once more Barbashin-Krasovsky's theorem to conclude
that the equilibrium point $P_{1}$ is asymptotically stable
and the set $\DD$ is contained inside its basin of attraction.

\begin{figure}[htbp]
\begin{centering}
\ins{360pt}{-94pt}{$x$}
\ins{214pt}{-10pt}{$y$}
\ins{170pt}{-100pt}{$-c_{0}$}
\ins{316pt}{-94pt}{$2c_{0}$}
\ins{306pt}{-50pt}{$\Gamma_{0}$}
\includegraphics*[angle=0,width=4in]{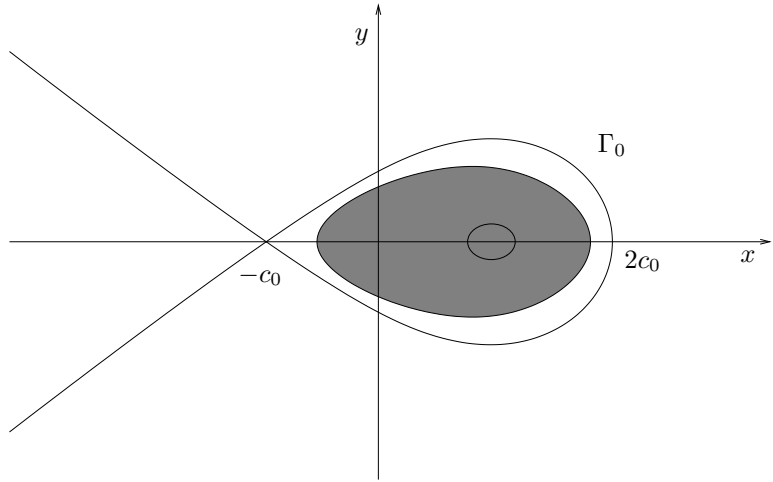}
\caption{\label{fig5-1}
Estimate of the basin of attraction (grey set) of the limit
cycle (solid line inside the grey set) for large $\g$. The curve
$\Gamma_{0}$ is the separatrix of the system described by
the differential equations $\ddot x + x^{2} = c_{0}^{2}$.}
\end{centering}
\end{figure}

If we go back to the original coordinates, we see that,
for $\g$ large enough, there is a quasi-periodic solution
$x_{0}(t)$ which moves very close to the point $(c_{0},0)$
(within a distance of order $\eps=1/\g$);
cf. figure~\ref{fig5-1}. Its basin of attraction contains
a large set whose boundary is at a distance
of order $\eps^{\be}$ from the separatrix $\Gamma_{0}$
of the system, that is from the curve of equation
$y = \pm \sqrt{2(2c_{0}^{3}-x^{3}+3c_{0}^{2}x)/3}$.

\zerarcounters
\section{Improved estimates for the basins of attraction}
\label{sec:6}

We start by constructing an invariant set $\AAA$,
valid for all $p\in\NNN$, which has the property that its area
is $O(\g)$. The construction is surely not optimal but
is included because it is, amongst many sets constructed,
one which grows with $\g$ in the vertical direction.
Constructing invariant sets containing the attracting orbit
with $\g$ considered to be large is not difficult, but finding
ones that grow with $\g$ was found to be less trivial. The
construction follows very closely the ideas exploited in~\cite{BDGM}.

Rewrite~(\ref{eq:1.4}), with $g(x)=x^{2p}$ as
\begin{equation}
\begin{cases} \dot x = y , \\
\dot y = f(\oo t) - \g y - x^{2p} , \end{cases}
\label{eq:6.1} \end{equation}
so that the vector field generated by the differential equation is
defined by $\vec{\phi}(t) = (y, f(\oo t) - \g y - x^{2p})$.
Let $f^{2p}\leq f(\oo t)\leq F^{2p}$. Note that in this
and the following section, $f$ is just a constant and not
the function $t\to f(\oo t)$. We adopt this notation
here to conform to that used in the analysis in~\cite{BDGM}.
For the same reason we denote in boldface vectors in $\RRR^{2}$,
so that $\vec{a}\cdot\vec{b}$ denotes the scalar product in $\RRR^{2}$.

The two vector fields, $\vec{\phi}_F = (y,F^{2p}-\g y-x^{2p})$
and $\vec{\phi}_f = (y, f^{2p} - \g y - x^{2p})$, have no explicit
time-dependence and also have the property that for all $t$, 
$\vec{\phi}(t) = \mu \vec{\phi}_F + (1-\mu)\vec{\phi}_f$ where $\mu$ is a
(time-dependent) scalar.

\begin{figure}[htbp]
\begin{centering}
\includegraphics*[angle=0,width=2.9in]{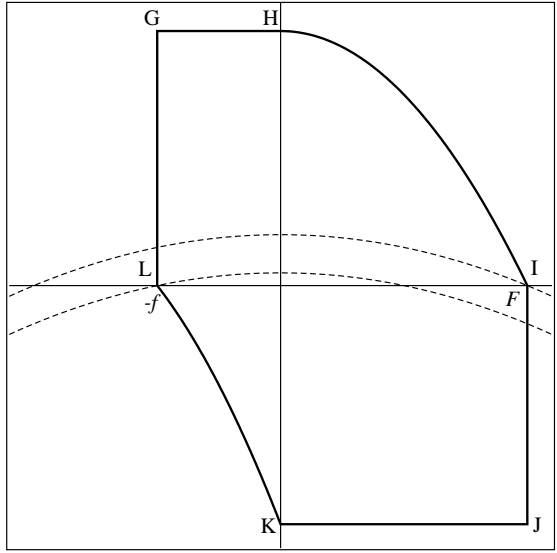}
\caption{\label{bigA}
The absorbing set GHIJKL, $\AAA$.
Also shown as dotted lines are the curves 
$P_F: y = (F^{2p} - x^{2p})/\g$ and $P_f: y = (f^{2p} - x^{2p})/\g$.}
\end{centering}
\end{figure}

Following~\cite{BDGM}, we let the boundary of invariant set $\AAA$
be a hexagon GHIJKL whose edges are straight lines, 
except for HI and KL. GH and JK are horizontal and LG and IJ are
vertical --- see figure~\ref{bigA}. The co-ordinates of points H and K
are $(0, y_H)$ and $(0, y_K)$ respectively. The dotted curves in
figure~\ref{bigA} are $P_F: y = (F^{2p} - x^{2p})/\g$ (upper)
and $P_f: y = (f^{2p} - x^{2p})/\g$ (lower). As shown in~\cite{BDGM},
only in the region between these curves is the sign of $\dot{y}$ ambiguous;
above $P_F$, $\dot{y} < 0$ and below $P_f$, $\dot{y} > 0$. Provided that
$y_H \geq F^{2p}/\g$, which will turn out to be automatically satisfied,
$\vec{\phi}(t)$ will be into $\AAA$ along GH, and point J being
below $P_f$ guarantees that $\vec{\phi}(t)$ will always be into KJ,
both results holding for all $t$. Furthermore, the sign of $\dot{x} = y$
guarantees that $\vec{\phi}(t)$ will be into LG and IJ, again for all $t$.
It therefore remains to prove that $\vec{\phi}(t)$ is into the
curved sides HI and KL for all $t$.

Let us define HI by $y = \lambda_1\left(F^{2p} - x^{2p}\right)$ 
with $\lambda_1\in\RRR$ positive; then the inward-pointing
normal, $\vec{n}_1 = (-2p\lambda_1 x^{2p-1}, -1)$. In order to prove that
$\vec{\phi}(t)$ is into HI for all time, we need only show that
$\vec{n}_1\cdot\vec{\phi}_F \geq 0$ for $x\in[0, F]$. We have
$\vec{n}_1\cdot \vec{\phi}_F = \left[F^{2p} - x^{2p}\right]
\left[\lambda_1(\g - 2p\lambda_1 x^{2p-1}) - 1\right]$ and,
since $F^{2p} - x^{2p}\geq 0$, and the second bracket in the scalar
product reaches its minimum value over $[0, F]$ at $x = F$, we require
\begin{equation}
2pF^{2p-1} \lambda_1^2 - \gamma\lambda_1 + 1\leq 0 ,
\label{eq:6.2} \end{equation}
in order for the scalar product to be non-negative.
Solving this quadratic and choosing the larger solution gives
\begin{equation}
\lambda_1 = \frac{\gamma}{4pF^{2p-1}}\left[1 + \sqrt{1 -
\frac{8pF^{2p-1}}{\gamma^2}}\right]
\label{eq:6.3} \end{equation}
provided $\lambda_1\in\RRR$. This is true for $\g^2\geq 8pF^{2p-1}$.
We now show that this condition on $\g$ also forces 
$y_H = \lambda_1 F^{2p} \geq F^{2p}/\g$, or $\lambda_1\geq 1/\g$,
to be satisfied. Letting $q = 4pF^{2p-1}$, the condition for $\lambda_1$
to be real is $\g^2/q\geq 2$; using this in the definition
of $\lambda_1$ gives $(\g^2/q)\left(1+\sqrt{1-2q/\g^2}\right)\geq 
2\left(1+\sqrt{1-2q/\g^2}\right)\geq 1$, which is clearly true.
Hence, boundary HI has been constructed in such a way that $\vec{\phi}(t)$
is into it for all time, provided that H is above $P_F$: this is true if 
$\lambda_1 \geq 1/\g$.  Note that $\lambda_1\sim \g/(2pF^{2p-1})$ as 
$\g\rightarrow\infty$, and so $y_H = \lambda_1 F^{2p} = O(\g)$.

We now define KL by $y = \lambda_2\left(f^{2p} - (x+2f)^{2p}\right)$
where $\lambda_2\in\RRR > 0$. This has inward 
normal $\vec{n}_2 = \left(2p\lambda_2(x+2f)^{2p-1}, 1\right)$
and the inequality we need to consider is now 
\begin{equation}
\vec{\phi}_f\cdot\vec{n}_2 = \lambda_2\left[(x+2f)^{2p} - f^{2p}\right]
\left[\g - 2p\lambda_2(x+2f)^{2p-1}\right] + f^{2p} - x^{2p}\geq 0
\label{eq:6.4} \end{equation}
for $x\in[-f, x_K]$ where $x_K$, satisfying $F\geq x_K> -f$,
is to be defined. The first term consists of a product of two terms,
the first of which is positive for $x\in(-f, \infty)$
and the second of which is a monotonically decreasing function of $x$ for
$x\in\mathbb R\setminus\{-2f\}$; it is zero at $x = x_0 =
(\g/2p\lambda_2)^{1/(2p-1)} - 2f$. Choosing $x_K = x_0 = 0$ ensures that the
product term in~(\ref{eq:6.4}) is non-negative and gives
\begin{equation}
\lambda_2 = \frac{\g}{2p(2f)^{2p-1}}.
\label{eq:6.5} \end{equation}
The last term in inequality~(\ref{eq:6.3}) is non-negative
for $x\in[-f, f]$ and so the scalar product is non-negative
for $\lambda_2$ as given above, provided that $x_K = 0$.
We also require point J to be below $p_f$. The $y$-co-ordinate
of J = $y_K = \lambda_2\left(f^{2p} - (2f)^{2p}\right)$ and
so this condition becomes 
\begin{equation}
\gamma^2 \geq \frac{p\left(F^{2p} - f^{2p}\right)}{f\left(1 - 4^{-p}\right)}.
\label{eq:6.6} \end{equation}

We can now state our result.

\begin{theorem}\label{thm:5}
Define $\AAA$ as the hexagonal set whose vertices
are G, H, I, J, K and L, where
\begin{eqnarray}
& & \qquad \mbox{G} = (-f, \lambda_1 F^{2p}),\;\;\;
\mbox{H} = (0, \lambda_1 F^{2p}), \;\;\;
\mbox{I} = (F, 0), \nonumber \\
& & \mbox{J} = (F, -\lambda_2 f^{2p}(4^p-1)), \;\;\;
\mbox{K} = (0, -\lambda_2 f^{2p}(4^p-1)), \;\;\;
\mbox{L} = (-f, 0) ,
\label{eq:6.7} \end{eqnarray}
and where $\lambda_1$ and $\lambda_2$ are given by
equations~~(\ref{eq:6.3}) and~~(\ref{eq:6.4}), respectively.
Let the edges LG, GH, IJ and JK
of $\AAA$ be straight lines and let HI be given
by $y = \lambda_1\left(F^{2p} - x^{2p}\right)$ and KL,
by $y = \lambda_2\left(f^{2p} - (x+2f)^{2p}\right)$.
Then, provided that 
\begin{equation}
\g^2 \geq \max\left(8pF^{2p-1}, \frac{p\left(F^{2p} -
f^{2p}\right)}{f\left(1 - 4^{-p}\right)}\right),
\label{eq:6.8} \end{equation}
set $\AAA$ is an invariant set containing the
limit cycle described by the solution $x_{0}(t)$ in the plane.
\end{theorem}

\begin{figure}[htbp]
\begin{centering}
\includegraphics*[angle=0,width=6.3in]{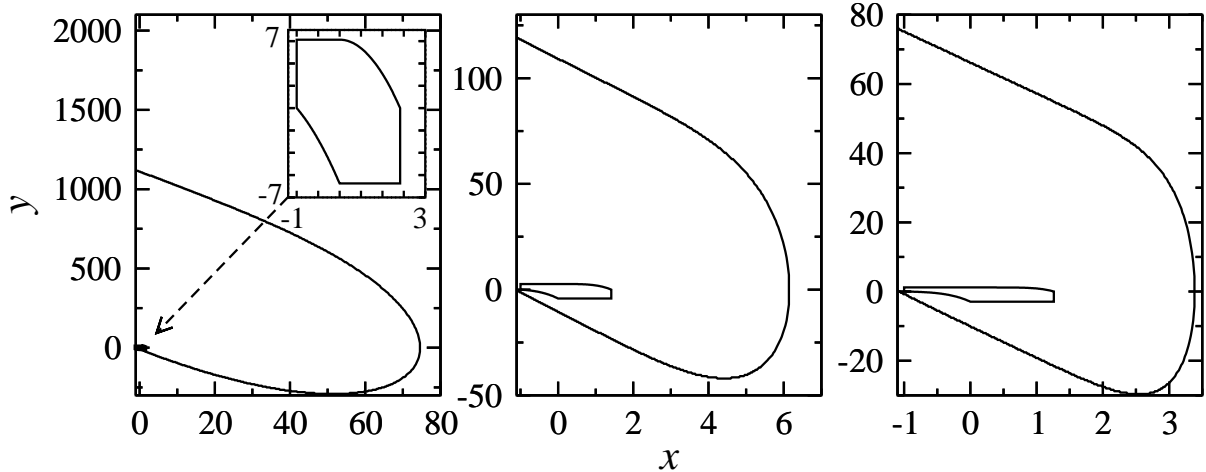}
\caption{\label{comp}
The set GHIJKL compared to the set $\mathcal B^c$
(cf.~\cite{BDGM}) with $1\leq f(t) = (5 + 3\sin t)/2 \leq 4$, $\gamma = 9$
and $p=1$ (left), $p=2$ (middle) and $p=3$ (right).
The set $\mathcal B$ is defined such that all
trajectories starting outside $\mathcal B$ are (numerically)
found to remain outside, hence it provides an estimate
for the basin of attraction of the periodic solution, $x_0(t)$.}
\end{centering}
\end{figure}

Figure~\ref{comp} compares the set GHIJKL with the 
actual basin of attraction for $p=1, 2$ and 3.

The existence of the set $\AAA$, with the properties
stated in Theorem~\ref{thm:5} allows us to improve the estimate
$\DD$ of the basin of attraction found in Section~\ref{sec:5}
for the varactor equation studied in~\cite{BDGM}, that is for
the system~(\ref{eq:7.1}), where $f(\oo t)=\al + \beta \sin t$,
with $\al>0$ and $|\beta|<\al$.
Define $\DD_{0} = \DD \cup \AAA$. First note that both
$\DD$ and $\AAA$ are strictly contained in $\DD_{0}$,
as follows easily from the remark that $\DD_{0}$ is inside
the separatrix $\Gamma_{0}$ (cf. the end of Section~\ref{sec:5}),
and at a distance $O(1/\sqrt{\g})$ from it.
The inclusion $\DD \subset \DD_{0}$ is obvious as $\AAA$
grows linearly at $\g$ in the vertical direction.
The inclusion $\AAA \not\subseteq \DD_{0}$ follows from the fact
that the separatrix intersects the $x$-axis in
$x=-c_{0}$ and $x=2c_{0}$, with $c_{0}=\sqrt{\al}$;
hence $x_{L}=-\sqrt{\al-\beta} > - c_{0}$ and
$x_{I}=\sqrt{\al+\beta} \le \sqrt{2\al} < 2 c_{0}$.
Then the set $\DD_{1} := \AAA \setminus (\AAA \cap \DD)$
is non-empty. Furthermore a trajectory starting in $\DD_{1}$
can go out of $\DD_{1}$ only by entering $\DD_{0}$ (by the invariance
of $\AAA$). This means that all trajectories starting in $\DD_{0}$
are attracted by the limit cycle described by $x_{0}(t)$.
 
\zerarcounters
\section{Blow up in finite time}
\label{sec:7}

We prove finite-time blow up for the system~(\ref{eq:6.1})
by first finding an invariant set ${\mathcal J}$
(Lemma~\ref{lem:7}) and then constructing an invariant subset,
${\mathcal S}(-X_0)$, of ${\mathcal J}$, the latter construction
being carried out in Lemma~\ref{lem:8}.
Within ${\mathcal S}(-X_0)$, a differential inequality must
hold and we show in Theorem~\ref{thm:6} that all
solutions of this differential inequality must blow up in finite time.

\begin{lemma} \label{lem:7}
Define $h(x) = 2px^{2p-1}\left(F^{2p} - x^{2p}\right) -
\g^2\left(F^{2p} - f^{2p}\right)$ and the set
\begin{equation}
{\mathcal{D}} = \{(x,y)| x\leq-\xi, (F^{2p} - x^{2p})/\g
\leq y \leq 0\} ,
\label{eq:7.1} \end{equation}
where $-\xi$ is the root of $h(x)$ with $-\xi < -F$. Then $\xi$ as
defined exists uniquely and ${\mathcal{D}}$ is an invariant set.
\end{lemma}

\prova The set ${\mathcal{D}}$ is clearly absorbing along its
horizontal, $y = 0$, $x\leq-\xi$, and vertical, $x = -\xi$,
$(F^{2p}- \xi^{2p})/\gamma \leq y\leq 0$ boundaries.
We therefore only need to prove that it is also absorbing along the
curved boundary, $y = (F^{2p} - x^{2p})/\gamma$ for $x\leq-\xi$.
The appropriate normal here is $\vec{n} = \left(2px^{2p-1}/\gamma, 1\right)$
and so we require $\vec{n}\cdot\vec{\phi}_f \geq 0$, which is equivalent
to proving that $h(x)\geq 0$, for $x\leq -\xi$.  

To this end, first note some elementary properties of $h(x)$:
(i) $h(-F) = \gamma^2\left(F^{2p} - f^{2p}\right) < 0$;
(ii) $h(x)$ has exactly three stationary points for $x\in\RRR$:
$h'(x) = 0$ for $x = 0$ and $x = x_\pm = \pm F[(2p-1)/(4p-1)]^{1/2p}$;
(iii) $h''(x_-) > 0$; and
(iv) $h(x)\rightarrow\infty$ as $x\rightarrow -\infty$.
From these, it becomes clear that $-F < x_- < 0$ and $h(x_-)$
is a minimum; and, in the light of (i) and (iv), $h(x)$ has exactly
one real root, $-\xi,\in(-\infty, -F)$. Hence, $h(x)\geq 0$
for $x\leq-\xi$ and the invariance of set ${\mathcal{D}}$ is proved.$\EP$

We now define, for $x\leq -X_0$, a curve 
$G(x, y): y = -b\left(-X_0 - x\right)^\rho$ with $-X_0\leq-\xi< 0$, 
$\rho = 3/2$, and $b>0$ to be found.  In order that
${\mathcal S}(-X_0)\subset {\mathcal J}$, $-X_0\leq -\xi$ and $b>0$.
In the proof of Theorem~\ref{thm:6}, it will be required
that ${\mathcal S}(-X_0)$ has infinite area for all $p\in\NNN$,
and hence $\rho<2$ --- otherwise, the curves $(F^{2p}- \xi^{2p})/\g$
and $G$ could intersect at some finite $x < -\xi$ when $p=1$.
In order that solutions blow up in {\em finite} time, $\rho>1$
will also be required: hence, we choose $\rho = 3/2$.

The proof of the invariance of ${\mathcal S}$ now follows.

\begin{lemma} \label{lem:8}
The set
\begin{equation}
{\mathcal S}(-X_0) = \{(x, y)| x \leq -\xi,\; U\geq y
\geq (F^{2p} - x^{2p})/\gamma\},
\label{eq:7.2} \end{equation}
where
\begin{equation}
U = \left\{\begin{array}{cll} -b\left(-X_0 - x\right)^{3/2} &
\mbox{\hskip 0.2in}& x \leq -X_0\\
0 &\mbox{\hskip 0.2in}& \mbox{otherwise}
\end{array}\right.
\label{eq:7.3} \end{equation}
with $-X_0\leq -\xi$ and $b$ sufficiently small, is an
invariant subset of $\mathcal J$, and contains points $(x, y)$
with $y\rightarrow-\infty$.
\end{lemma}

\prova The fact that the vector field is into all the boundaries of
$\mathcal S$ except $y = -b\left(-X_0 - x\right)^{3/2}$ has been proved in
Lemma~\ref{lem:7}; therefore, we only need to consider 
the boundary $y = -b\left(-X_0 - x\right)^{3/2}$.
The appropriate normal is that which points into ${\mathcal S}$,
this being $\vec{n} = (3b/2)\sqrt{-X_0 - x}, -1)$.
The correct choice for the vector field here is
$\vec{\phi}_F = (y, F^{2p} - \gamma y - x^{2p})$, which gives
\begin{equation}
\vec{n}\cdot\vec{\phi}_F = -b(-X_0 - x)^{3/2}\left[(3b/2)\sqrt{-X_0 -
x}+\g\right] + x^{2p} - F^{2p} \geq 0,\;\;x\leq-X_0.
\label{eq:7.4} \end{equation}
Since $-X_0 - x \geq 0$, we can substitute $v^2 = -X_0 - x$, giving
\begin{equation}
L(v) =  (X_0+v^2)^{2p} - F^{2p} -bv^3(3bv/2 + \g) \geq 0, \;\; v\geq 0.
\label{eq:7.5} \end{equation}
Let $l(x) = x^p$; then the mean value theorem states that,
for $b > a$, $l(b) - l(a) = (b-a)p c^{p-1}$, where $c\in(a, b)$.
Applying this to the first two terms in equation~~(\ref{eq:7.5}),
we have
\begin{eqnarray}
L(v) & \geq & \left[(X_0+v^2)^2 - F^2\right]
pF^{2p-2} - bv^3(3bv/2 + \g) \nonumber \\
& = & v^4(pF^{2p-2} - 3b^2/2) - \g b v^3 +
2X_0 pF^{2p-2}v^2 +(X_0^2 - F^2)pF^{2p-2} \nonumber \\
& \geq & v^2\left[ v^2(pF^{2p-2} - 3b^2/2) - \g b v +
2X_0 pF^{2p-2}\right] = v^2 M(v) ,
\label{eq:7.6} \end{eqnarray}
where we have used $X_0 > F$ to obtain the last inequality. Hence,
we need to show that for $b>0$ sufficiently small, $M(v)\geq0$
for $v\geq 0$. First, write $M(v) = a_2 v^2 - a_1b v + a_0$, $a_1, a_0 >0$, 
and let $b$ be sufficiently small that $a_2>0$.
Let $\tilde{M}(v) = M(v)/a_2 = v^2 - 2c_1 bv + c_0$,
with $c_1, c_0 > 0$. Then $\tilde{M}(v) = (v - b c_1)^2 + c_0 -
b^2 c_1^2$, and it is plain that $b$ can be chosen to be
small enough that $\tilde{M}(v)$, and so $L(v)$, are non-negative
for $v\geq 0$.

We will also require $\mathcal S(-X_0)$ to extend to infinite
negative $y$ values; it has this property provided that the
curves $y = -b\left(-X_0 - x\right)^{3/2}$ and
$y = (F^{2p} - x^{2p})/\g$ nowhere intersect for $y\leq 0$.
This, too, is clearly true if $b$ is small enough.$\EP$

It is possible that some of the conditions applied in the above proof
could be relaxed, but a `better' invariant set ${\mathcal S}$ is not
required in  the proof of the following theorem.

\begin{theorem}\label{thm:6}
For all $\vec{x}_0 = (x_0, y_0)\in {\mathcal J}$, 
$\exists -X_0\leq-\xi$ which is such that $\vec{x}_0\in {\mathcal S}(-X_0)$.
All solutions starting from such an $\vec{x}_0$ blow up in finite time.
\end{theorem}

\prova Let $\vec{x}_0 = \left(x(0), y(0)\right) = (x_0, y_0)
\in {\mathcal J}$. Then there always exists an $-X_0 \leq -\xi$
such that $\vec{x}_0$ is in a subset ${\mathcal S}(-X_0)$
of ${\mathcal J}$: choose any $-X_0 > x_0$. Also, since $\dot{x} \leq 0$ in
${\mathcal J}$, $x\leq x_0$ and so $-X_0 - x \geq -X_0 - x_0 > 0$.
Additionally, since $\vec{x}_0\in {\mathcal S}(-X_0)$, the
differential inequality $-b(-X_0 - x)^{3/2} \geq y = \dot{x}
\geq (F^{2p} - x^{2p})/\gamma$ applies for all
time $t \geq 0$, by the invariance of ${\mathcal S}(-X_0)$.

In fact, only the upper bound is important here, and with the
substitution $u = -X_0 - x > 0$ and $u_0 = -X_0 -x_0 > 0$,
this becomes $\dot{u} \geq b u^{3/2}$. Integrating gives
\begin{equation}
\int_{u_0}^u {\rm d} \psi \, \psi^{-3/2} \geq b\int_0^t {\rm d} s ,
\label{eq:7.7} \end{equation}
giving $u_0^{-1/2} - u^{-1/2}\geq bt/2$, which, after re-arranging, gives
\begin{equation}
u\geq \left[ u_0^{-1/2} - bt/2\right]^{-2}\;\;\;
\mbox{ for } t\in[0, t_\infty) ,
\label{eq:7.8} \end{equation}
where $b t_\infty = 2 u_0^{-1/2}$. The above inequality
shows that $x(t)$ tends to $-\infty$ within a finite time
$t \leq t_\infty$. By the invariance of ${\mathcal S}(-X_0)$, $x(t)$
cannot tend to $-\infty$ without also $y(t) \rightarrow -\infty$,
and so finite time blow up is proven.$\EP$

Note that for $p=1$ there is a set $\mathcal B$ (cf.
figure~\ref{comp}) such that all solutions starting from $\mathcal B$
reach ${\mathcal J}$ in a finite time: this was proved in~\cite{BDGM}.
Therefore for $p=1$ Theorem \ref{thm:6} shows that all
trajectories starting from $\mathcal B$ blow up in finite time.

\zerarcounters
\section{Conclusions, extensions and open problems}
\label{sec:8}

We conclude this paper with a list of open problems
(some already mentioned in the previous sections).

The first one concerns possible extensions of the
proof of Theorem~\ref{thm:3} to the case of more general
polynomials of the form~~(\ref{eq:4.1}). A natural question
is under which conditions there is still a global attractor,
in these cases, when the dissipation coefficient is large enough.

A characterisation of the set $\BB$ can be given in some concrete cases,
such as that of the varactor equation considered in Section~\ref{sec:5}.
Its diameter is of order $1$. We have seen in Section~\ref{sec:6}
that we can improve the estimate by obtaining a set whose size
increases linearly in $\g$ in the vertical direction,
but in such a way that it is still expected to be strictly
included inside the actual basin of attraction.
It would be worthwhile to attempt constructions of 
sets contained inside the basins of attraction
that are as large as possible.

We also leave as an open problem for the varactor equation
the proof that any bounded solution is attracted
by $x_{0}(t)$. On the basis of numerical simulations,
we conjecture that this is the case.

Another interesting problem is whether one can weaken the hypotheses
on the function $g$, both for determining the existence of a
quasi-periodic solution with the same frequency vector as the forcing and,
in that case, for proving its uniqueness and attractivity.

Finally, extensions to higher dimensional cases would be desirable.

\appendix

\zerarcounters
\section{Proof of~(\ref{eq:3.10})}
\label{app:A}

One has $\dot P = \partial_{\xi} P \, \dot \x + \partial_{t}P$
and $\dot Q = \partial_{\xi} Q \, \dot \x$, so that
\begin{equation}
\frac{\dot P}{P} - \frac{Q}{Q} =
\left( \frac{\partial_{\xi}P}{P} -
\frac{\partial_{\xi}Q}{Q} \right) \dot \xi +
\frac{\partial_{t}P}{P} = \frac{\partial_{\xi}P - \partial_{\xi}Q}{P} +
\frac{\partial_{\xi}Q}{PQ} \left( Q - P \right) +
\frac{\partial_{t}P}{P} .
\label{eq:A.1} \end{equation}
One can write
\begin{eqnarray}
\partial_{t}P & = & \sum_{j=1}^{2p}
\left( \begin{matrix} 2p+1 \\ j \end{matrix} \right)
j \xi^{2p+1-j} x_{0}^{j-1}(t)\, \dot x_{0}(t) , \\ \nonumber
Q - P & = & \sum_{j=1}^{2p}
\left( \begin{matrix} 2p+1 \\ j \end{matrix} \right)
\xi^{2p+1-j} \left( \al - x_{0}(t) \right)^{j} , \\ \nonumber
\partial_{\xi}P - \partial_{\xi} Q & = & \sum_{j=1}^{2p}
\left( \begin{matrix} 2p+1 \\ j \end{matrix} \right)
\left( 2p+1-j \right)  \xi^{2p-j}
\left( x_{0}(t) - \al \right)^{j} ,
 \label{eq:A.2} \end{eqnarray}
where $\dot x_{0}(t)=O(1/\g)$ and $x_{0}(t)-\al = O(1/\g)$.

Finally $|x_{0}(t)| \le 2|\al|$ for all $t\in\RRR$ if $\g$ is
large enough, and both $\xi^{2p+1-j}/P$ and $\xi^{2p+1-j}/Q$
tend to zero as $\xi\to\infty$ for $j\geq 1$. Hence~(\ref{eq:3.10})
follows, with the constants $B_{1}$ and $B_{2}$ depending on $p$
but not on $\g$.

\zerarcounters
\section{Initial data in I and III}
\label{app:B}

Take an initial datum $z=(\xi,y)$ in I. If $y=0$ then $\dot\xi =0$
and $\dot y = - \xi F(\xi,x_{0}(t))<0$, so that the trajectory
enters II. If $y>0$ then $\dot \xi >0$ and $\dot y<0$.

Moreover $\partial_{\xi} (\xi F(\xi,x))=
\partial_{\xi}(x+\xi)^{2p+1}=(2p+1)\,(x+\xi)^{2p} \geq 0$
for all $x\in\RRR$, so that, by using the fact that $\xi(t) \ge \xi(0)$
as long as $(\xi(t),y(t))$ remains in I, one has in I
\begin{equation}
\xi F(\xi,x_{0}(t)) \ge \inf_{t\in\RRR}
\xi F(\xi,x_{0}(t)) \ge \xi(0) \inf_{t\in\RRR}
F(\xi(0),x_{0}(t)) \geq c > 0 ,
\label{eq:B.1} \end{equation}
where we used that $F(\xi,x)$ is strictly greater than a positive
constant for $x \neq 0$ (see the proof of Lemma~\ref{lem:1}).
Therefore we obtain
\begin{equation}
\dot y \leq - \g y - c ,
\label{eq:B.2} \end{equation}
which implies that $y(t)$ reaches the $\xi$-axis in a finite time.

Analogously one discusses the case of initial data $z$ in III.

\zerarcounters
\section{On the curves $\boldsymbol \CC_{1}$ and $\boldsymbol \CC_{2}$}
\label{app:C}

Call $\TT$ the subset of II outside $\SSS$.

Define $\CC_{1}$ as a continuous curve in II such that in $\TT$
it is given by the graph of the function $\xi \to -\xi^{2p+1}/4\g$.
In~(\ref{eq:3.3}) one can write $\dot y = \g (-y+g(\xi,t))$,
with $g(\xi,t)$ defined in Lemma~\ref{lem:3}.
By Lemma~\ref{lem:3}, in $\TT$ one has $g(\xi,t) \leq 
-\xi^{2p+1}/2\g$,
so that at all points in $\TT$ above $\CC_{1}$ one has
$$ -y+g(\xi,t)) = |y| + g(\xi,t) \leq \frac{1}{4\g}\,\xi^{2p+1} -
\frac{1}{2\g}\,\xi^{2p+1} \leq - \frac{1}{4\g}\,\xi^{2p+1} , $$
hence $\dot y < 0$.

Define $\CC_{2}$ as a continuous curve in II such that in $\TT$
it is given by the graph of the function $\xi \to - 4 \xi^{2p+1}/\g$.
By Lemma~\ref{lem:3}, one has $g(\xi,t) \geq -2\xi^{2p+1}/\g$,
so that in all points of $\TT$ below $\CC_{2}$ one has
$y \leq - 4 \xi^{2p+1}/\g \leq 2 g(\xi,t)$, hence
$$ -\g y \geq \g(- y + g(\xi,t)) \geq - \g y /2 , $$
so that $ -\g y \geq \dot y \geq - \g y/2$.
In terms of the rescaled variables $(X,Y)$ this yields
$Y' \equiv \Psi(X,Y)$, with
\begin{equation}
-\g\eps^{p} Y \geq \Psi(X,Y) \ge -\g\eps^{p} Y/2 ,
\label{eq:C.1} \end{equation}
as asserted after~(\ref{eq:3.21}).

The point $P$ is given by the intersection of the curve $Y_{1}(X)=-
4 X^{2p+1}/\g\eps^{p}$ withe the level curve $\Gamma$. Hence
$$ \frac{1}{2} = \frac{1}{2} \left( \frac{4X_{P}^{2p+1}}{\g\eps^{p}}
\right)^{2} + X\,F(X,\eps x_{0}(0)) =
\frac{2^{4p+1}}{\g^{2}\eps^{2p}} X_{P}^{4p+2} +
X_{P}^{2p+1} + O(\eps\, X^{2p}) , $$
hence $X_{P}=O((\g\eps^{p})^{1/(2p+1)})$.

Now consider the solution of~(\ref{eq:3.19}) with initial datum
$Z(T_{2})=(X(T_{2}),Y(T_{2}))$. We want to check that the solution
remains below $\CC_{2}$ until it crosses the $Y$-axis.
The solution of
\begin{equation}
\begin{cases} X' = Y , \\
Y' = \Psi(X,Y) , \end{cases}
\label{eq:C.2} \end{equation}
with $\Psi(X,Y)$ satisfying the bounds~(\ref{eq:C.1}),
moves below the line with slope $-\g\eps^{p}$
passing through $Z(T_{2})$, that is below the line of equation
\begin{equation}
Y_{1}(X) = Y_{0} - \g \eps^{p} X
\label{eq:C.3} \end{equation}
with $Y_{0}$ determined by the request that for $X=\overline X \equiv
X(T_{2})$ one has
\begin{equation}
Y_{0} - \g \eps^{p} \overline X = - \frac{4}{\g \eps^{p}}
\overline X^{2p+1} ,
\label{eq:C.4} \end{equation}
where the graph of $-4X^{2p+1}/\g \eps^{p}$ describes the curve 
$\CC_{2}$ in the coordinates $(X,Y)$. By using that $\overline X$
is close to $X_{P}$ one realises that $Y_{0}$ in~(\ref{eq:C.4})
has to be negative. In turn this implies that the line
of equation~(\ref{eq:C.3}) is below the curve $\CC_{2}$,
so that also the assertion after~(\ref{eq:3.23}) is proved.

\zerarcounters
\section{Variations in finite times for quasi-integrable systems}
\label{app:D}

The system obtained from~(\ref{eq:3.19}) by replacing $\Psi(X,Y)$
with $\Psi_{1}(X,Y)$ is an integrable Hamiltonian system,
with Hamiltonian~(\ref{eq:3.21}). For $\eps=0$ the Hamiltonian
reduces to
\begin{equation}
H_{0}(X,Y) = \frac{1}{2} Y^{2} + \frac{1}{2p+2} X^{2p+2} ,
\label{eq:D.1} \end{equation}
which can be written in terms of the action-angle variables
$(I,\varphi)$ as
\begin{equation}
H_{0}(X,Y) = \HHH_{0}(I) = c_{p} I^{(2n+2)/(n+2)} ,
\label{eq:D.2} \end{equation}
where $c_{p}$ is a suitable $p$-dependent positive constant.
By taking into account the other terms of the vector field, we obtain
\begin{equation}
H(X,Y) = \HHH(I) = c_{p} I^{(2n+2)/(n+2)} + O(I^{(2n+1)/(n+2)}) .
\label{eq:D.3} \end{equation}
The equations obtained by adding to $\Psi_{1}(X,Y)$ the vector
field $\Psi_{2}(X,Y)$ are still Hamiltonian, and are described
by the non-autonomous Hamiltonian
\begin{equation}
\HHH(I) + \HHH_{1}(I,\varphi,t) ,
\label{eq:D.4} \end{equation}
with $\HHH_{0}$ given as in~(\ref{eq:D.3}) and $\HHH_{1}$ of order
$\eps^{p+1}$ as long as the action variables remain of order $1$.

The corresponding equations of motion are
$$ \begin{cases}
\dot I = - \partial_{\varphi} \HHH_{1}(I,\varphi) , \\
\dot \varphi = \omega_{0}(I) + \partial_{I} \HHH_{1}(I,\varphi) ,
\end{cases} $$
with $\omega_{0}(I)=\partial_{I}\HHH_{0}(I)$. Then one
immediately realises that in a time of order $1$ the action
variables remain close to their initial values.
In turn this implies that also the angle variables are changed
by order $\eps^{p+1}$ with respect their unperturbed values.
In terms of the original coordinates $(X,Y)$ this means that
the solution remains within a distance $O(\eps^{p+1})$ with respect
the unperturbed value.


\end{document}